\documentclass[11pt,a4paper]{article} 
\usepackage[margin=25mm]{geometry}
\usepackage{cite,amsmath,amsfonts,amssymb,stmaryrd,latexsym, amsthm, cases, comment, mathrsfs, tikz,  enumerate, mathtools, multicol, xcolor, xspace}
\usepackage{enumitem}
\usepackage[margin=10pt,font=small,labelfont=bf, labelsep=period]{caption}

\usepackage[colorlinks]{hyperref}

\newcommand\nc\newcommand
\nc\rnc\renewcommand

\usetikzlibrary{decorations.markings}
\usetikzlibrary{arrows,matrix}
\usepgflibrary{arrows}
  
\tikzset{->-/.style={decoration={
  markings,
  mark=at position #1 with {\arrow{{latex}}}},postaction={decorate}}}

\nc\bit{\begin{itemize}}
\nc\eit{\end{itemize}}
\nc\ben{\begin{enumerate}[label=\textup{(\roman*)},leftmargin=7mm]}
\nc\BEN{\begin{enumerate}[label=\textup{(\Roman*)},leftmargin=7mm]}
\nc\bena{\begin{enumerate}[label=\textup{(\alph*)},leftmargin=7mm]}
\nc\een{\end{enumerate}}
\nc\bmc{\begin{multicols}}
\nc\emc{\end{multicols}}
\rnc\implies{\ \Rightarrow\ }
\rnc\iff{\ \Leftrightarrow\ }

\nc\Sing{\operatorname{Sing}}
\nc\supp{\operatorname{supp}}
\nc\id{\operatorname{id}}
\nc\N{\mathbb N}
\nc\ba{{\bf a}}
\nc\bb{{\bf b}}
\nc\bone{{\bf 1}}
\nc\bn{{\bf n}}
\nc\bm{{\bf m}}
\nc\set[2]{\{#1:#2\}}
\nc\bigset[2]{\big\{#1:#2\big\}}
\nc\AND{\qquad\text{and}\qquad}
\nc\ANd{\quad\text{and}\quad}
\nc\ANDSIM{\qquad\text{and similarly}\qquad}
\nc\COMMA{,\qquad}
\nc\COMMa{,\quad}
\nc\WHERE{\qquad\text{where}\qquad}
\nc\io\iota
\nc\si\sigma
\nc\Om\Omega
\nc\Ga\Gamma
\nc\De\Delta
\nc\sub\subseteq
\nc\sm\setminus
\nc\mt\mapsto
\nc\D{\mathcal D}
\nc\M{\mathcal M}
\nc\tsharp{_{\scriptscriptstyle{\op}}^\sharp}
\rnc\S{\mathcal S}
\nc\ol\overline
\nc\ul\underline
\nc\wh\widehat
\nc\ulwh[1]{\ul{\hspace{0.7truemm}\wh{#1}\hspace{0.7truemm}}}
\nc\comp{^\complement}
\nc{\ldb}{[\hspace{-0.5truemm}[}
\nc{\rdb}{]\hspace{-0.5truemm}]}
\nc\com{^\complement}
\nc\es\varnothing

\nc\fTC{^{\hspace{-0.5truemm}\oast}}

\nc\la\langle
\nc\ra\rangle
\nc{\pres}[2]{\la {#1} : {#2} \ra} 
\nc{\Spres}[2]{\textup{\textsf{Sgp}}\la {#1} : {#2} \ra} 
\nc{\Mpres}[2]{\textup{\textsf{Mon}}\la {#1} : {#2} \ra} 
\nc{\Cpres}[2]{\textup{\textsf{Cat}}\la {#1} : {#2} \ra} 
\nc{\TCpres}[2]{\textup{\textsf{TCat}}\la {#1} : {#2} \ra} 
\nc\trans[1]{\left(\begin{smallmatrix}#1\end{smallmatrix}\right)}

\let\oldproofname=\proofname
\rnc{\proofname}{\rm\bf{\oldproofname}}
\nc\pf{\begin{proof}}
\nc\epf{\end{proof}}
\nc\epfres{\hfill\qed}
\nc\epfreseq{\tag*{\qed}}

\nc\XS{X_{\Sing(\I_n)}}
\nc\RS{R_{\Sing(\I_n)}}
\nc\pS{\phi_{\Sing(\I_n)}}
\nc\XSM{X_{\Sing(M_0^n)}}
\nc\RSM{R_{\Sing(M_0^n)}}
\nc\pSM{\phi_{\Sing(M_0^n)}}
\nc\XMS{X_{M\wr\Sing(\I_n)}}
\nc\RMS{R_{M\wr\Sing(\I_n)}}
\nc\pMS{\phi_{M\wr\Sing(\I_n)}}

\nc{\pfitem}[1]{\medskip\noindent #1.}
\nc{\firstpfitem}[1]{#1.}
\nc{\pfcase}[1]{\medskip\noindent {\bf Case #1.}}
\newcommand{\uv}[1]{\fill (#1,2)circle(.17);}
\newcommand{\lv}[1]{\fill (#1,0)circle(.17);}
\newcommand{\uvx}[2]{\fill (#1,2)circle(#2);}
\newcommand{\lvx}[2]{\fill (#1,0)circle(#2);}
\newcommand{\uvs}[1]{{\foreach \x in {#1} { \uv{\x}}}}
\newcommand{\lvs}[1]{{\foreach \x in {#1} { \lv{\x}}}}
\newcommand{\uvxs}[2]{{\foreach \x in {#1} { \uvx{\x}{#2}}}}
\newcommand{\lvxs}[2]{{\foreach \x in {#1} { \lvx{\x}{#2}}}}
\newcommand{\stline}[2]{\draw(#1,2)--(#2,0);}
\newcommand{\stlines}[1]{{\foreach \x/\y in {#1} { \stline{\x}{\y} }}}
\newcommand{\uvc}[2]{\fill[#2] (#1,2)circle(.17);}
\newcommand{\lvc}[2]{\fill[#2] (#1,0)circle(.17);}
\newcommand{\uvcs}[2]{{\foreach \x in {#1} { \uvc{\x}{#2}}}}
\newcommand{\lvcs}[2]{{\foreach \x in {#1} { \lvc{\x}{#2}}}}
\newcommand{\stlinec}[3]{\draw[#3](#1,2)--(#2,0);}
\newcommand{\stlinecs}[2]{{\foreach \x/\y in {#2} { \stlinec{\x}{\y}{#1} }}}
\nc\udotted[2]{\draw[dotted](#1+.5,2)--(#2-.5,2);}
\nc\ddotted[2]{\draw[dotted](#1+.5,0)--(#2-.5,0);}
\nc{\custpartn}[3]{{\lower1.4 ex\hbox{
\begin{tikzpicture}[scale=.3]
\uvs{#1}
\lvs{#2}
#3 \end{tikzpicture}
}}}

\numberwithin{equation}{section}

\newtheorem{theorem}[equation]{Theorem}
\newtheorem{lemma}[equation]{Lemma}

\theoremstyle{definition}

\newtheorem{remark}[equation]{Remark}

\newtheorem{eg}[equation]{Example}

\DeclareMathOperator{\im}{im}
\DeclareMathOperator{\dom}{dom}

\newcommand{\br}{\mathbf{r}}
\newcommand{\bd}{\mathbf{d}}
\newcommand{\lam}{\lambda}

\newcommand{\invmon}{\mathcal{I}_n}
\newcommand{\n}{\mathbf{n}}
\newcommand{\op}{\oplus}
\newcommand{\C}{\mathcal{C}}

\newcommand{\I}{\mathcal{I}}
\newcommand{\al}{\alpha}
\newcommand{\be}{\beta}

\newcommand{\U}{\rotatebox[origin=c]{180}{$U$}}

\begin{document}

\title{\vspace{-1cm}Presentations for wreath products involving symmetric inverse monoids and categories}

\author{Chad Clark and James East\footnote{Supported by ARC Future Fellowship FT190100632.}\\[3mm]
{\it\small Centre for Research in Mathematics and Data Science,}\\
{\it\small Western Sydney University, Locked Bag 1797, Penrith NSW 2751, Australia.}\\[3mm]
{\tt\small Chad.Clark\,@\,WesternSydney.edu.au}, {\tt\small J.East\,@\,WesternSydney.edu.au}}

\date{}

\maketitle

\begin{abstract}
Wreath products involving symmetric inverse monoids/semigroups/categories arise in many areas of algebra and science, and presentations by generators and relations are crucial tools in such studies. The current paper finds such presentations for $M\wr\I_n$, $M\wr\Sing(\I_n)$ and $M\wr\I$.  Here $M$ is an arbitrary monoid, $\I_n$ is the symmetric inverse monoid, $\Sing(\I_n)$ its singular ideal, and $\I$ is the symmetric inverse category.



\end{abstract}

\tableofcontents

\section{Introduction}

Semigroups whose elements are fundamental mathematical entities such as partitions, transformations and relations have a natural representation using diagrams. Since the work of Brauer in the~1930s~\cite{brauer1937algebras}, so-called diagram semigroups have appeared in a number diverse settings, including topology, invariant theory, quantum mechanics, representation theory, statistical mechanics and many other branches of mathematics and science.  See for example \cite{HR2005,Kauffman1987, francis2021brauer, lehrer2012second,temperley1971relations, kauffman2007mathematics,jones1994potts, martin1994temperley, Martin2008, phyloalg2}, and especially \cite{Martin2008} for an overview/survey; more examples and references can be found in the introductions to \cite{east2020presentations,EG2017}.  

Finding presentations of diagram/transformation semigroups and their singular parts, which are ideals formed by complementing the group of units,  has been of considerable interest since at least the late 1950s.  This began with the work of A\u{\i}zen\v{s}tat \cite{aizenshtat1962defining,aizenshtat1958defining} and Popova \cite{popova1961defining,popova1962defining}, and can be traced back to Moore's presentation for the symmetric group as a Coxeter group \cite{moore}; cf.~\cite{Humphreys1990}.  With the notable exceptions of the semigroups of (partitioned) binary relations \cite{plemmons1970semigroup,martin2013partitioned}, presentations are known for a great many important diagram semigroups. These include the full/partial transformation monoids, (dual) symmetric inverse monoids, Temperley-Lieb/Kauffman monoids, (partial) Brauer monoids, Motzkin monoids, (rook) partition monoids, and the singular ideals of many of the above; see for example \cite{BDP2002,aizenshtat1962defining,aizenshtat1958defining,popova1961defining,popova1962defining, East2015,kudryavtseva2006presentations, maltcev2006presentation, East2021,East2019,East2018,HLP2013,east2020presentations,FitzGerald2003}.

The symmetric inverse monoid $\invmon$ is the set of all partial permutations of the set ${\bn=\{1,\ldots,n\}}$, under the operation of relational composition.  This monoid plays the same role in the theory of inverse semigroups as the symmetric group $\S_n$ plays in group theory; for more details on this, see for example \cite{Lawson1998,lipscomb1996symmetric}. The monoid~$\invmon$ is also known in the literature as the rook monoid because of the alternate characterisation of partial permutations by $\{0,1\}$-matrices with at most one non-zero entry in each row and each column; such matrices are in one-to-one correspondence with placements of non-attacking rooks on an $n\times n$ chess board.  While the representation theory of the rook monoid is well studied (see for example \cite{solomon2002representations, Munn1964}), there has also been substantial recent interest in the representation theory of wreath products~$G\wr\I_n$, where~$G$ is a group (definitions are given below).  See especially the work of Steinberg \cite{steinberg2008mobius, steinberg2016representation} and Mazorchuk and Srivastava~\cite{mazorchuk2021jucys}.  The paper \cite{mazorchuk2021jucys} concerns the important special case that $G=C_r$ is a finite cyclic group, which stems from analogous wreath products $C_r\wr\S_n$ and their deformations, the cyclotomic Hecke algebras \cite{DJ1998}.  See also \cite{Brookes2021}, which computes the congruence lattice of $G\wr\I_n$ for an arbitrary group $G$.  As the title of that paper suggests, $G\wr\I_n$ is in fact isomorphic to the monoid of partial automorphisms of a free $G$-act of rank $n$, providing another source of motivation for studying such products.

The current article studies more general structures still, and concerns \emph{transformational} wreath products $M\wr\I_n$ for an arbitrary \emph{monoid} $M$, as studied for example in \cite{Brookes2021, CDEGZ} (full definitions are given in Subsection \ref{subsect:wreath}). We also consider:
\bit
\item $M\wr\Sing(\I_n)$, where $\Sing(\I_n)$ is the singular ideal of $\I_n$, consisting of all \emph{strictly} partial permutations, and
\item $M\wr\I$, where $\I$ is the category of all partial bijections $\bm\to\bn$ where $m$ and $n$ range over $\N = \{0, 1, \dots\}$.
\eit
Our main results give presentations by generators and relations for each of these structures. Although several general results exist on presentations of wreath products of monoids (see \cite{dombi2009generators,robertson2002finite,howie1994constructions, robertson2003finite}), these are not applicable to $M\wr\I_n$ as they concern different kinds of wreath products to the transformational types considered here.  As for~$M\wr\I$, we are not aware of any previous results on presentations for wreath products involving categories, but we will be aided by tools developed in \cite{east2020presentations}, which allow presentations for certain classes of categories to be constructed from presentations of their endomorphism monoids.  The situation for $M\wr\Sing(\I_n)$ is rather more complex, as $\Sing(\I_n)$ is not a monoid, and it is well known that semigroups tend to behave quite badly even for much simpler constructions such as Cartesian/direct products \cite{RRW1998}.  As a case in point, an entire 41-page paper \cite{feng2019presentations} has been devoted to the related wreath product $M\wr\Sing(\mathcal T_n)$, where $\mathcal T_n$ is the full transformation monoid (consisting of all self-maps of $\n$).  While the recent paper \cite{CDEGZ} deals with a vast class of semigroup products, and proves many general results, examples such as $M\wr\Sing(\I_n)$ were singled out as a special kind of difficult case, and none of the results proved in \cite{CDEGZ} apply here.  Thus, we hope that these methods of working with non-monoid wreath products will be of use in the study of other such complex structures.

The article is organised as follows. We begin with some preliminary material in Section~\ref{sec:preliminaries}, which includes the definition of the wreath products that are our focus.  Section \ref{sect:MIn} gives presentations for the monoid $M\wr\I_n$; see Theorems \ref{thm:MIn} and \ref{thm:MIn2}, each of which extends a known presentation for~$\I_n$~\cite{popova1961defining,Gilbert2006}.  In Section \ref{sec:category}, we apply the results of Section \ref{sect:MIn} and the general machinery developed in \cite{east2020presentations} to obtain a category presentation (Theorem~\ref{thm:MI}) of $M \wr \I$, which is then used to obtain a tensor presentation (Theorem~\ref{thm:MI2}) for $M\wr\I$. This tensor presentation involves a particularly compact set of generators and relations, which is afforded by the tensor (a.k.a monoidal) structure of $M \wr \I$. This structure allows morphisms between large objects to be built from morphisms between smaller objects, as made precise below. Finally, we treat the singular wreath product $M\wr\Sing(\I_n)$ in Section \ref{sec:wreathprodsing}. This involves the most technically-demanding arguments, and an intermediate result of independent interest gives a presentation for the (singular) semigroup $M_0^n\sm M^n$, where here $M_0$ is the semigroup obtained by adjoining a new zero element to the monoid~$M$; see Theorem \ref{thm:M}.  This is then combined with the presentation for $\Sing(\I_n)$ from~\cite{East2015} to prove our final main result, Theorem \ref{thm:MSIn}, which gives a presentation for $M\wr\Sing(\I_n)$. 

\subsection*{Acknowledgement}

The authors would like to thank the anonymous referee for their careful reading of the original manuscript and their valuable suggestions.

\section{Preliminaries}\label{sec:preliminaries}

We require only the most basic notions from semigroup/monoid theory, as may be found in monographs such as \cite{Howie1995,Higgins1992,CP1961}.  We now gather the preliminary material and background results we need concerning categories (Subsection \ref{subsect:cat}), symmetric inverse monoids and categories (Subsection~\ref{subsect:I}), wreath products (Subsection \ref{subsect:wreath}) and presentations (Subsection \ref{subsect:pres}).

\subsection{Categories}\label{subsect:cat}

As in \cite{east2020presentations}, we are interested in very special kinds of categories.  Throughout, the term \emph{category} will always mean a small, strict tensor (a.k.a.~monoidal) category with object set $\N$; these terms will all be defined below.  Any such category $\C$ will be identified with its morphism set.  The domain and range of $a\in\C$ are denoted by $\bd(a)$ and $\br(a)$, respectively.  Morphisms are composed left-to-right, so that $ab=a\circ b$ is defined if and only if $\br(a)=\bd(b)$, in which case $\bd(ab)=\bd(a)$ and $\br(ab)=\br(b)$.  For $m,n\in\N$ we write $\C_{m,n}=\set{a\in\C}{\bd(a)=m,\ \br(a)=n}$ for the set of all morphisms $m\to n$.  For $n\in\N$ we write $\C_n=\C_{n,n}$ for the endomorphism monoid at the object $n$, and we denote the identity of this monoid by $\io_n$ (or a similar symbol).  To say that $\C$ is a \emph{strict tensor category} means that it has an additional (totally defined) operation $\op:\C\times\C\to\C$ for which 
\bit
\item $(\C,\op)$ is a monoid with identity $\io_0$,
\item $\bd(a\op b)=\bd(a)+\bd(b)$ and $\br(a\op b)=\br(a)+\br(b)$ for all $a,b\in\C$,
\item $\io_m\op\io_n=\io_{m+n}$ for all $m,n\in\N$,
\item $(a\circ b)\op(c\circ d) = (a\op c)\circ(b\op d)$ for all $a,b,c,d\in\C$ with $\br(a)=\bd(b)$ and $\br(c)=\bd(d)$.
\eit
A \emph{congruence} on a category $\C$ (as above) is an equivalence relation $\si$ on $\C$ for which
\bit
\item $(a,b)\in\si \implies \bd(a)=\bd(b)$ and $\br(a)=\br(b)$ for all $a,b\in\C$,
\item $(a,b)\in\si \implies (xa,xb),(ay,by)\in\si$ for all $a,b,x,y\in\C$ whenever the stated products are defined.
\eit
The \emph{quotient category} $\C/\si$ is the set of all $\si$-classes under the induced operation.  We say $\si$ is a \emph{tensor congruence} if it also preserves $\op$, meaning that 
\bit
\item $(a,b)\in\si \implies (x\op a,x\op b),(a\op y,b\op y)\in\si$ for all $a,b,x,y\in\C$.
\eit
For a set $\Om\sub\C\times\C$ satisfying $(a,b)\in\Om\implies\bd(a)=\bd(b)$ and $\br(a)=\br(b)$, we denote by $\Om^\sharp$ the congruence on $\C$ generated by $\Om$.  So $\Om^\sharp$ is the smallest congruence on $\C$ containing $\Om$.  We also denote by $\Om\tsharp$ the tensor congruence generated by $\Om$; note that $\Om^\sharp\sub\Om\tsharp$, but the inclusion is sometimes strict.

Given categories $\C$ and $\D$ (as above), a \emph{category morphism} $\phi:\C\to\D$ is a functor that acts as the identity on objects, meaning that $\bd(a\phi)=\bd(a)$ and $\br(a\phi)=\br(a)$ for all $a\in\C$.  The \emph{kernel} of~$\phi$ is the congruence
\[
\ker(\phi) = \bigset{(a,b)\in\C\times\C}{a\phi=b\phi},
\]
and we have $\C/{\ker(\phi)} \cong\im(\phi)$.

\subsection{Symmetric inverse monoids and categories}\label{subsect:I}

A particularly important category for us is the \emph{symmetric inverse category} $\I$, defined as follows.  For $n\in\N$ let $\bn=\{1,\ldots,n\}$, interpreting this to be empty when $n=0$.  For $m,n\in\N$ we denote by $\I_{m,n}$ the set of all injective partial maps $\bm\to\bn$: i.e., all injective functions $A\to\bn$ with $A\sub\bm$.  We then define
\[
\I = \bigset{(m,\al,n)}{m,n\in\N,\ \al\in\I_{m,n}}.
\]
For $m,n,k\in\N$, and for $\al\in\I_{m,n}$ and $\be\in\I_{n,k}$, we define
\[
\bd(m,\al,n)=m \COMMA \br(m,\al,n)=n \AND (m,\al,n)\circ(n,\be,k) = (m,\al\be,k),
\]
where $\al\be=\al\circ\be\in\I_{m,k}$ is the ordinary relational composition.  To avoid clutter in our notation, we will typically identify an element $(m,\al,n)$ of $\I$ with the partial map $\al\in\I_{m,n}$ itself, but regard~$m$ and $n$ as `encoded' in~$\al$, writing $\bd(\al)=m$ and $\br(\al)=n$.  In this way, the morphism sets of $\I$ are the $\I_{m,n}$ ($m,n\in\N$), and the endomorphism monoids are the \emph{symmetric inverse monoids} $\I_n=\I_{n,n}$.  

A partial bijection $\al\in\I_{m,n}$ is typically represented by a graph consisting of a row of upper vertices labelled $1,\ldots,m$, a row of lower vertices labelled $1,\ldots,n$, and an edge from upper vertex~$i$ to lower vertex $j$ if $i\in\dom(\al)$ and $i\al=j$.  Here as usual $\dom(\al)$ denotes the domain of~$\al$; we also write $\im(\al)$ for the image of $\al$.  Figure \ref{fig:I} gives representations of the partial bijections
\[
\al=\trans{1&2&3&4&5&6\\-&-&4&1&-&7}\in\I_{6,8} \AND \be=\trans{1&2&3&4&5&6&7&8\\2&-&1&-&5&-&4&-}\in\I_{8,7},
\]
as well as their composition $\al\be=\trans{1&2&3&4&5&6\\-&-&-&2&-&4}\in\I_{6,7}$, noting that an entry of `$-$' in position $i$ indicates that $i$ does not belong to the domain of the mapping.  This diagrammatic interpretation of composition in $\I$ will help with many of the calculations to follow.  We assume, unless otherwise specified, that the vertices in each row are labelled in left-to-right ascending order and generally omit these labels for convenience.

\begin{figure}[ht]
\begin{center}
\begin{tikzpicture}[scale=.6]
\begin{scope}[shift={(0,0)}]	
\uvxs{1,...,6}{.15}
\lvxs{1,...,8}{.15}
\stlines{3/4,4/1,6/7}
\draw(0.6,1)node[left]{$\al=$};
\draw[-{latex}](10.5,-1)--(12.5,-1);
\foreach \x in {1,...,8} {\draw[dotted] (\x,0)--(\x,-2);}
\end{scope}
\begin{scope}[shift={(0,-4)}]	
\uvxs{1,...,8}{.15}
\lvxs{1,...,7}{.15}
\stlines{1/2,3/1,5/5,7/4}
\draw(0.6,1)node[left]{$\be=$};
\end{scope}
\begin{scope}[shift={(14,-2)}]	
\uvxs{1,...,6}{.15}
\lvxs{1,...,7}{.15}
\stlines{4/2,6/4}
\draw(7.4,1)node[right]{$=\al\be$};
\end{scope}
\end{tikzpicture}
\caption{Partial bijections $\al\in\I_{6,8}$ and $\be\in\I_{8,7}$, and their composition $\al\be\in\I_{6,7}$.}
\label{fig:I}
\end{center}
\end{figure}

The category $\I$ also has a tensor operation.  For $\al\in\I_{m,n}$ and $\be\in\I_{k,l}$, we define the mapping ${\al\op\be\in\I_{m+k,n+l}}$ by
\[
x(\al\op\be) = \begin{cases}
x\al &\text{if $1\leq x\leq m$ and $x\in\dom(\al)$,}\\
(x-m)\be+n &\text{if $m+1\leq x\leq m+k$ and $x-m\in\dom(\be)$,}\\
- &\text{otherwise.}
\end{cases}
\]
Diagrammatically, $\al\op\be$ is obtained by stacking $\al$ and $\be$ horizontally, as in Figure \ref{fig:I2}.

\begin{figure}[ht]
\begin{center}
\begin{tikzpicture}[scale=.52]
\begin{scope}[shift={(0,0)}]	
\uvcs{1,...,6}{red}
\lvcs{1,...,8}{red}
\stlinecs{red}{3/4,4/1,6/7}
\draw(0.6,1)node[left]{${\color{red}\al}=$};
\draw[->](9.5,-1)--(11.5,-1);
\end{scope}
\begin{scope}[shift={(0,-4)}]	
\uvcs{1,...,7}{blue}
\lvcs{1,...,6}{blue}
\stlinecs{blue}{1/2,3/1,5/5,7/4}
\draw(0.6,1)node[left]{${\color{blue}\be}=$};
\end{scope}
\begin{scope}[shift={(12,-2)}]	
\uvcs{1,...,6}{red}
\lvcs{1,...,8}{red}
\stlinecs{red}{3/4,4/1,6/7}
\uvcs{7,...,13}{blue}
\lvcs{9,...,14}{blue}
\stlinecs{blue}{7/10,9/9,11/13,13/12}
\draw(14.4,1)node[right]{$={\color{red}\al}\op{\color{blue}\be}$};
\end{scope}
\end{tikzpicture}
\caption{Partial bijections $\al\in\I_{6,8}$ and $\be\in\I_{7,6}$, and their sum $\al\op\be\in\I_{13,14}$.}
\label{fig:I2}
\end{center}
\end{figure}

\begin{comment}
It is worth noting that $\I$ is an inverse category in the sense of \cite{Kastl1979} and \cite[Section 2.3.2]{CL2002}.  The inverse of $\al\in\I_{m,n}$ is its ordinary inverse map $\al^{-1}\in\I_{n,m}$, and we have $\al=\al\al^{-1}\al$ and $\al^{-1}=\al^{-1}\al\al^{-1}$ for all $\al\in\I$.  The category $\I$ has various other natural structural properties; for example, it is a PROP in the sense of \cite[Section 24]{MacLane1965}.  These additional structures will not play an explicit role here, but \cite[Remark 4.22]{east2020presentations} explains how the PROP relations can be incorporated into presentations for $\I$.
\end{comment}

We have already noted that the endomorphism monoids of the category $\I$ are the symmetric inverse monoids $\I_n=\I_{n,n}$.  The group of units of $\I_n$ is the \emph{symmetric group} $\S_n$.  Because the monoids here are finite, the singular part $\Sing(\I_n)=\I_n\sm\S_n$ is itself a semigroup; indeed, it is an ideal of $\I_n$.  The set $\S=\bigcup_{n\in\N}\S_n$ of all units of $\I$ is of course a subcategory, but $\Sing(\I)=\I\sm\S$ is not a category; for example, if $\al$ and $\be$ denote the unique elements of $\I_{0,1}$ and $\I_{1,0}$, respectively, then $\al,\be\in\Sing(\I)$, while $\al\be=\es$ is the unique element of $\I_0=\S_0$.

\subsection{Wreath products}\label{subsect:wreath}

In this subsection we define the transformational wreath products that will be the objects of our study. 

Let $M$ be a monoid with identity $1$ and let $0$ be a symbol not belonging to $M$.  We denote by $M_0 = M\cup\{0\}$ the monoid obtained from $M$ by adjoining $0$ as a multiplicative zero element, even if $M$ already had a multiplicative zero element.  For $n\in\N$ we write $M_0^n$ for the set of all $n$-tuples over $M_0$, which is a monoid with identity $\bone_n=(1,\ldots,1)$.  Note that $M_0^0=\{\es\}$ consists only of the empty tuple.  We also write
\[
\M = M_0^0 \cup M_0^1\cup M_0^2\cup\cdots
\]
for the set of all finite tuples over $M_0$, which is trivially a category (in our sense).  For $\ba\in \M$ we write $\bd(\ba)=\br(\ba)$ for the unique $n\in\N$ such that $\ba\in M_0^n$.  For $\ba=(a_1,\ldots,a_m)\in M_0^m$ and $\bb=(b_1,\ldots,b_n)\in M_0^n$ we define $\ba\op\bb=(a_1,\ldots,a_m,b_1,\ldots,b_n)\in M_0^{m+n}$.

The \emph{support} of $\ba=(a_1,\ldots,a_n)\in M_0^n$ is defined by
\[
\supp(\ba) = \set{i\in\bn}{a_i\not=0}.
\]
Since $\supp(\ba\cdot\bb) = \supp(\ba)\cap\supp(\bb)$ for all $\ba,\bb\in M_0^n$, it follows that $\supp:M_0^n\to2^{\bn}$ is in fact a monoid (sur)morphism, where $2^\bn$ denotes the power set $\set{A}{A\sub\bn}$, considered as a semilattice (semigroup of commuting idempotents) under $\cap$.

Let $\al\in\I_{m,n}$ and let $\ba=(a_1,\ldots,a_n)\in M_0^n$.  We define
\[
{}^\al{\ba} = (b_1,\ldots,b_m)\in M_0^m \WHERE b_i = \begin{cases}
a_{i\al} &\text{if $i\in\dom(\al)$}\\
0 &\text{otherwise.}
\end{cases}
\]
Note that ${}^\al\ba$ is defined if and only if $\bd(\ba)=\br(\al)$, in which case $\bd({}^\al\ba)=\bd(\al)$.  
Diagrammatically, one can represent a tuple $\ba=(a_1,\ldots,a_n)\in M_0^n$ as a row of $n$ vertices, with the $i$th vertex labelled~$a_i$.  The action $\ba\mt{}^\al\ba$ can then be calculated by `sliding' the relevant entries of $\ba$ up the edges of $\al$; Figure \ref{fig:act} gives an example with $\al=\trans{1&2&3&4&5&6\\-&-&4&1&-&7}\in\I_{6,8}$.

\begin{figure}[h]
\begin{center}
\scalebox{.85}{
\begin{tikzpicture}[scale=.95]
\tikzstyle{vertex}=[circle,draw=black, fill=white, inner sep = 0.07cm]
\draw[-{latex}](9,1)--(11,1);
\begin{scope}[shift={(0,0)}]	
\uvxs{1,...,6}{.1}
\lvxs{1,...,8}{.1}
\foreach \x/\y in {3/4,4/1,6/7} {\draw[->-=0.5] (\y,0)--(\x,2);}
\foreach \x in {1,...,8} {\draw[dotted] (\x,0)--(\x,-1);}
\end{scope}
\begin{scope}[shift={(0,-3)}]	
\foreach \x in {1,...,8} {\node[vertex] () at (\x,2){$a_{\x}$};}
\node[left] () at (0,2) {$\ba$};
\draw[-{latex}](0,2)--(.5,2);
\end{scope}
\begin{scope}[shift={(11,0)}]	
\uvxs{1,...,6}{.1}
\lvxs{1,...,8}{.1}
\stlines{3/4,4/1,6/7}
\foreach \x in {1,...,6} {\draw[dotted] (\x,2)--(\x,3);}
\end{scope}
\begin{scope}[shift={(11,1)}]	
\foreach \x/\y in {3/4,4/1,6/7} {\node[vertex] () at (\x,2){$a_{\y}$};}
\foreach \x in {1,2,5} {\node[vertex,minimum size=6.3mm] () at (\x,2){$0$};}
\node[right] () at (7,2) {${}^\al\ba$};
\draw[-{latex}](7,2)--(6.5,2);
\end{scope}
\end{tikzpicture}
}
\caption{The action of $\al\in\I_{6,8}$ on $\ba\in M_0^8$.}
\label{fig:act}
\end{center}
\end{figure}

One can easily check that
\[
{}^\al(\ba\cdot\bb) = {}^\al\ba\cdot{}^\al\bb \AND {}^{\al\be}\ba = {}^\al({}^\be\ba)
\]
for all $\al,\be\in\I$ and $\ba,\bb\in \M$ for which the above expressions are defined.  We also have ${}^{\id_n}\ba=\ba$ for all $\ba\in M_0^n$.  However, we do not necessarily have ${}^\al\bone_n=\bone_m$ for $\al\in\I_{m,n}$; rather, ${}^\al\bone_n$ is the $m$-tuple over $\{0,1\}$ with support equal to $\dom(\al)$.

\begin{figure}[h]
\begin{center}
\scalebox{.75}{
\begin{tikzpicture}[scale=.9]
\tikzstyle{vertex}=[circle,draw=black, fill=white, inner sep = 0.06cm]
\begin{scope}[shift={(0,0)}]	
\uvxs{1,...,6}{.12}
\lvxs{1,...,8}{.12}
\stlines{3/4,4/1,6/7}
\draw(0.6,1)node[left]{$(\ba,\al)=$};
\draw[-{latex}](9.5,-1)--(11.5,-1);
\foreach \x in {1,...,8} {\draw[dotted] (\x,0)--(\x,-2);}
\foreach \x in {3,4,6} {\node[vertex] () at (\x,2){$a_{\x}$};}
\end{scope}
\begin{scope}[shift={(0,-4)}]	
\uvxs{1,...,8}{.12}
\lvxs{1,...,7}{.12}
\stlines{1/2,3/1,5/5,7/4}
\draw(0.6,1)node[left]{$(\bb,\be)=$};
\foreach \x in {1,3,5,7} {\node[vertex] () at (\x,2){$b_{\x}$};}
\end{scope}
\begin{scope}[shift={(12,-2)}]	
\uvxs{1,...,6}{.12}
\lvxs{1,...,7}{.12}
\stlines{4/2,6/4}
\draw(7.4,1)node[right]{$=(\ba,\al)\circ(\bb,\be)$};
\foreach \x/\y in {4/1,6/7} {\node[vertex] () at (\x,2){$a_{\x}b_{\y}$};}
\end{scope}
\end{tikzpicture}
}
\caption{Elements of $M\wr\I$ (left) and their composition (right).}
\label{fig:I3}
\end{center}
\end{figure}
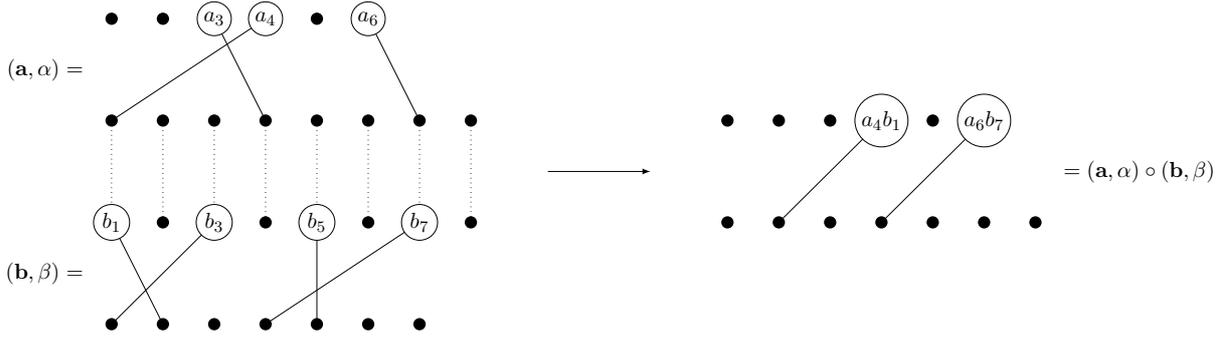

The \emph{wreath product} $M\wr\I$ is defined as follows.  As a set we have
\[
M\wr\I = \bigset{(\ba,\al)\in \M\times\I}{\bd(\ba)=\bd(\al),\ \supp(\ba)=\dom(\al)}.
\]
For $(\ba,\al)\in M\wr\I$ we define $\bd(\ba,\al)=\bd(\al)$ and $\br(\ba,\al)=\br(\al)$.  For $(\ba,\al),(\bb,\be)\in M\wr\I$ with $\br(\al)=\bd(\be)$ we define 
\[
(\ba,\al)\circ(\bb,\be) = (\ba\cdot{}^\al\bb,\al\be).
\]
This composition is well defined (i.e., $\bd(\ba\cdot{}^\al\bb)=\bd(\al\be)$ and ${\supp(\ba\cdot{}^\al\bb)=\dom(\al\be)}$ for such $\ba,\bb,\al,\be$) and associative, and $M\wr\I$ forms a category.  The tensor operation is given by
\[
(\ba,\al)\op(\bb,\be) = (\ba\op\bb,\al\op\be).
\]
Again there is a handy diagrammatic representation of the elements and operations of $M\wr\I$.  We represent $(\ba,\al)$ by picturing $\al$ in the way described above, and we label the $i$th upper vertex with $a_i$.  We always omit the label $0$ from an upper vertex, and sometimes also do this with the label $1$; one can always tell if an unspecified upper vertex label is $0$ or $1$ according to whether the corresponding element of $\bm$ belongs to $\dom(\al)$ or not.  Figure \ref{fig:I3} gives an example calculation of a composition.  The tensor operation again corresponds to a horizontal stacking, as in Figure \ref{fig:I2}.  Of course when $M=\{1\}$ is trivial, $M\wr\I\cong\I$.

The endomorphism monoids of $M\wr\I$ are the ordinary monoid wreath products
\[
M\wr\I_n = \bigset{(\ba,\al)\in M_0^n\times \I_n}{\supp(\ba)=\dom(\al)},
\]
as studied for example in \cite{Brookes2021,CDEGZ}.

Note that $M\wr\I$ contains natural isomorphic copies of $\M$ and $\I$, and we will identify these with subcategories of $M\wr\I$ as follows.  For $\ba\in M_0^n$ and $\al\in\I_{m,n}$ we identify
\[
\ba \equiv (\ba,\id_{\supp(\ba)}) \AND \al \equiv (\bone_{\dom(\al)},\al),
\]
where $\id_{\supp(\ba)}\in\I_{\bd(\ba)}$ is the identity map on $\supp(\ba)$ and where $\bone_{\dom(\al)}\in M_0^{\bd(\al)}$ denotes the $\bd(\al)$-tuple over $\{0,1\}$ with support equal to $\dom(\al)$.

\begin{remark}\label{rem:aabb}
With the above identifications, we note that 
\[
(\ba,\al) = \ba\cdot\al \qquad\text{for any $(\ba,\al)\in M\wr\I$.}
\]
More generally, if $\ba\in M_0^m$ and $\al\in\I_{m,n}$ are arbitrary, then $\ba\cdot\al = (\bb,\be)$, where:
\bit
\item $\bb$ is obtained from $\ba$ by replacing $a_i$ by $0$ for any $i\not\in\dom(\al)$, and
\item $\be$ is the restriction of $\al$ to $\supp(\ba)$.
\eit
\end{remark}

For any $n\in\N$ and any subsemigroup $S$ of $\I_n$ we also have the wreath product
\[
M\wr S = \bigset{(\ba,\al)\in M_0^n\times S}{\supp(\ba)=\dom(\al)},
\]
which is a subsemigroup of $M\wr\I_n$. Later we will be especially concerned with the case that $S$ is the singular subsemigroup $\Sing(\I_n)=\I_n\sm\S_n$.

\subsection{Presentations}\label{subsect:pres}

We now fix the notation we will be using for (monoid, semigroup, category and tensor category) presentations.  We also state two results from \cite{east2020presentations} that will be required.

Fix an alphabet $X$ and let $X^*$ be the \emph{free monoid} over $X$.  So $X^*$ consists of all words over~$X$, including the empty word, which we will typically denote by $\io$.  Let $R\sub X^*\times X^*$, and let $R^\sharp$ be the congruence on $X^*$ generated by $R$.  We say a monoid $M$ has \emph{presentation} $\Mpres XR$ if $M\cong X^*/R^\sharp$: i.e., if there exists a monoid surmorphism $X^*\to M$ with kernel $R^\sharp$.  If $\phi$ is such a surmorphism, we say $M$ has \emph{presentation $\Mpres XR$ via $\phi$}.  Elements of $X$ and $R$ are called \emph{generators} and \emph{relations}, respectively, and a relation $(u,v)\in R$ is typically displayed as an equation: $u=v$.  (The same conventions hold for the other kinds of presentations discussed below.)  We sometimes use the notation $\Mpres XR$ to denote the monoid $X^*/R^\sharp$ itself.  

\emph{Semigroup presentations} are defined analogously in terms of the \emph{free semigroup} ${X^+=X^*\sm\{\io\}}$.  Relations in a semigroup presentation $\Spres XR$ always involve non-empty words, meaning that ${R\sub X^+\times X^+}$.

For any alphabet $X$, we denote by $\ell=\ell_X:X^*\to\N$ the length function.  The following basic result will be useful on a number of occasions.

\begin{lemma}\label{lem:XY}
Let $X$ and $Y$ be disjoint alphabets, let $\sim$ be a congruence on $(X\cup Y)^*$ or $(X\cup Y)^+$, and suppose one of the following two conditions holds:
\ben
\item \label{XY1} for all $x\in X$ and $y\in Y$, $yx\sim uv$ for some $u\in X^*$ and $v\in Y^*$ with $\ell(u)\leq1$,
\item \label{XY2} for all $x\in X$ and $y\in Y$, $yx\sim uv$ for some $u\in X^*$ and $v\in Y^*$ with $\ell(v)\leq1$.
\een
Then for all $w\in (X\cup Y)^+$ we have $w\sim uv$ for some $u\in X^*$ and $v\in Y^*$.
\end{lemma}

\pf
By symmetry, we assume that \ref{XY1} holds, and we proceed by induction on $k=\ell(w)$.  The result being clear for $k=1$, we assume that $k\geq2$, so that $w=w'z$ for some $w'\in(X\cup Y)^+$ and $z\in X\cup Y$.  By induction, we have $w'\sim u'v'$ for some $u'\in X^*$ and $v'\in Y^*$, and so $w\sim u'v'z$.  If $z\in Y$ then we take $u=u'$ and $v=v'z$, so we now assume that $z\in X$.  Repeated application of condition \ref{XY1} gives $v'z\sim u''v$ for some $u''\in X^*$ and $v\in Y^*$ (with $\ell(u'')\leq1$), so $w\sim u'v'z\sim u'u''v$, and we take $u=u'u''$.
\epf

The assumption that $\ell(u)\leq1$ in condition \ref{XY1} of Lemma \ref{lem:XY} cannot be dropped (and similarly for condition \ref{XY2}).  For example, suppose $X=\{x\}$, $Y=\{y\}$, and that $\sim$ is the congruence on $\{x,y\}^*$ generated by the single relation $yx=xxyy$.  Then one can show by induction that any word $\sim$-equivalent to $yxx$ must have the form $x^{2k}\cdot yxx\cdot y^{2k}$ or $x^{2k}\cdot xxyyx\cdot y^{2k}$ for some $k\in\N$.  In particular, the conclusion of Lemma \ref{lem:XY} does not hold in this case.

For category presentations, we must work with digraphs and paths in place of alphabets and words.  Keeping in mind the conventions about categories from Subsection \ref{subsect:cat}, let $\Ga$ be a digraph with vertex set $\N$, possibly with multiple/parallel edges, and possibly with loops.  We identify $\Ga$ with its edge set, and denote the source and target of an edge $x\in\Ga$ by $\bd(x)$ and~$\br(x)$ respectively.  The \emph{free category} over $\Ga$ is the set~$\Ga^*$ of all paths in $\Ga$ under concatenation (where defined).  The empty path at $n\in\N$ will be denoted by~$\io_n$.  Every other path can be thought of as a word of the form $w=x_1\cdots x_k$, where $k\geq1$ and $x_1,\ldots,x_k\in\Ga$, and where $\br(x_i)=\bd(x_{i+1})$ for all $1\leq i<k$.  For such a word/path, we have $\bd(w)=\bd(x_1)$ and $\br(w)=\br(x_k)$.  Now let $\Om\sub\Ga^*\times\Ga^*$ be a set of pairs of paths, such that $\bd(u)=\bd(v)$ and $\br(u)=\br(v)$ for all $(u,v)\in\Om$.  We say a category~$\C$ (over~$\N$) has \emph{presentation} $\Cpres\Ga\Om$ if~$\C\cong\Ga^*/\Om^\sharp$: i.e., if there exists a surmorphism $\Ga^*\to \C$ with kernel~$\Om^\sharp$.  If $\phi$ is such a surmorphism, we say $\C$ has \emph{presentation $\Cpres \Ga\Om$ via $\phi$}.

Finally, we recall the formulation of (strict) tensor category presentations from \cite{east2020presentations}.  Let $\De$ be a digraph with vertex set $\N$, again identified with its edge set.  We will denote the \emph{free tensor category} over~$\De$ by~$\De\fTC$.  It consists of all terms constructed in the following way:
\begin{enumerate}[label=\textup{(T\arabic*)},leftmargin=12mm]
\item \label{T1} All empty paths $\io_n$ ($n\in\N$) are terms, with $\bd(\io_n)=\br(\io_n)=n$, acting as identities of $\De\fTC$.
\item \label{T2} All edges $x\in\De$ are terms, with $\bd(x)$ and $\br(x)$ the source and target of $x$, respectively.
\item \label{T3} If $s$ and $t$ are terms, and if $\br(s)=\bd(t)$, then the formal expression $s\circ t$ is a term, with $\bd(s\circ t)=\bd(s)$ and $\br(s\circ t)=\br(t)$.
\item \label{T4} If $s$ and $t$ are terms, then the formal expression $s\op t$ is a term, with $\bd(s\op t)=\bd(s)+\bd(t)$ and $\br(s\op t)=\br(s)+\br(t)$.
\een
Note that \ref{T1}--\ref{T4} describe the elements of $\De\fTC$, while \ref{T3} and \ref{T4} also give the definition of the $\circ$ and $\op$ operations.  Now let $\Xi\sub\De\fTC\times\De\fTC$ be a set of pairs of terms, such that $\bd(u)=\bd(v)$ and $\br(u)=\br(v)$ for all $(u,v)\in\Xi$.  We say a tensor category $\C$ (over $\N$) has presentation $\TCpres\De\Xi$ if~$\C\cong\De\fTC/\Xi\tsharp$: i.e., if there is a surmorphism $\De\fTC\to \C$ with kernel~$\Xi\tsharp$.  If $\phi$ is such a surmorphism, we say $\C$ has presentation $\TCpres\De\Xi$ via $\phi$.

As an example, we state the tensor presentation for $\I$ from \cite[Section 4.5]{east2020presentations}.  Let $\De_\I$ be the digraph with three edges denoted $X$, $U$ and $\U$, with
\[
\bd(X)=\br(X) = 2 \COMMA \bd(U)=\br(\U)=1 \AND \br(U)=\bd(\U)=0,
\]
and let $\Xi_\I$ be the following set of relations, where for convenience we write $I=\io_1\in\De\fTC$:
\begin{align}
\label{I1} &X\circ X=\io_2 \COMMA (X\op I)\circ(I\op X)\circ(X\op I) = (I\op X)\circ(X\op I)\circ(I\op X),\\
\label{I2} &\U\circ U = \io_0  \COMMA X\circ(U\op I) = I\op U \COMMA (\U\op I)\circ X = I\op \U  .
\end{align}

\begin{theorem}\label{thm:I}
The symmetric inverse category $\I$ has presentation $\TCpres{\De_\I}{\Xi_\I}$ via
\[
\epfreseq
X \mt \custpartn{1,2}{1,2}{\stline12\stline21} \COMMA
U \mt \custpartn{1}{}{\lvc{1}{white}} \COMMA
\U \mt \custpartn{}{1}{} .
\]
\end{theorem}

Note for example that $X\op I\mt \custpartn{1,2,3}{1,2,3}{\stline12\stline21\stline33}$ and $I\op X\mt  \custpartn{1,2,3}{1,2,3}{\stline11\stline23\stline32}$, so that the second relation in~\eqref{I1} is an instance of a so-called braid relation:
\[
\begin{tikzpicture}[scale=.3]
\uvs{1,2,3}
\lvs{1,2,3}
\stlines{1/2,2/1,3/3}
\begin{scope}[shift={(0,-2)}]
\uvs{1,2,3}
\lvs{1,2,3}
\stlines{1/1,2/3,3/2}
\node () at (5,1) {$=$};
\end{scope}
\begin{scope}[shift={(0,-4)}]
\uvs{1,2,3}
\lvs{1,2,3}
\stlines{1/2,2/1,3/3}
\end{scope}
\begin{scope}[shift={(6,0)}]
\uvs{1,2,3}
\lvs{1,2,3}
\stlines{1/1,2/3,3/2}
\end{scope}
\begin{scope}[shift={(6,-2)}]
\uvs{1,2,3}
\lvs{1,2,3}
\stlines{1/2,2/1,3/3}
\end{scope}
\begin{scope}[shift={(6,-4)}]
\uvs{1,2,3}
\lvs{1,2,3}
\stlines{1/1,2/3,3/2}
\end{scope}
\end{tikzpicture}
\]
The reader is invited to draw diagrams to check that the other relations listed in \eqref{I1} and \eqref{I2} hold.

The main general results from \cite{east2020presentations} provide mechanisms for building presentations for a category out of presentations for its endomorphism monoids.  Here we state special cases of these results that are sufficient for our purposes.  With this in mind, we fix a tensor category $\C$ (as above), and we denote the identities of $\C$ by $\ol\io_n\in\C_n$ ($n\in\N$).  In addition, we assume that all of the following conditions hold: 
\begin{enumerate}[label=\textup{(A\arabic*)}]
\item \label{A1} Every hom-set $\C_{m,n}$ ($m,n\in\N$) is non-empty.  (This actually follows from the next assumption, but we list it anyway, for ease of comparison with \cite{east2020presentations}.)
\item \label{A2} For every $n\in\N$, there exist $\ol\lam_n\in\C_{n,n+1}$ and $\ol\rho_n\in\C_{n+1,n}$ such that $\ol\lam_n\ol\rho_n=\ol\io_n$.
\item \label{A3} For every $n\in\N$, the endomorphism monoid $\C_n$ has presentation $\Mpres{X_n}{R_n}$ via ${\phi_n:X_n^*\to\C_n}$.  We assume the alphabets $X_n$ ($n\in\N$) are pairwise disjoint, and we denote by $\io_n$ the empty word over $X_n$.
\een
We now define a digraph $\Ga$ with vertex set $\N$ and edge set 
\[
\set{\lam_n,\rho_n}{n\in\N} \cup \bigcup_{n\in\N}X_n,
\]
with sources and targets given by
\[
\bd(x)=\br(x)=\bd(\lam_n)=\br(\rho_n)=n \ANd \br(\lam_n)=\bd(\rho_n)=n+1 \qquad\text{for all $n\in\N$ and $x\in X_n$.}
\]
We define a morphism $\phi:\Ga^*\to\C$ by
\[
\lam_n\phi=\ol\lam_n \COMMA \rho_n\phi=\ol\rho_n \AND x\phi=x\phi_n \qquad\text{for all $n\in\N$ and $x\in X_n$.}
\]
For $w\in\Ga^*$ we write $\ol w=w\phi\in\C$.  

We now assume that we have a set of relations $\Om\sub\Ga^*\times\Ga^*$ such that, writing ${\sim}$ for the congruence $\Om^\sharp$ on $\Ga^*$, the following all hold:
\begin{enumerate}[label=\textup{(A\arabic*)}]\addtocounter{enumi}{3}
\item \label{A4} For every relation $(u,v)\in\Om$, we have $\ol u=\ol v$.
\item \label{A5} For all $n\in\N$, $\Om$ contains the relation $\lam_n\rho_n=\io_n$, and a relation of the form $\rho_n\lam_n=w_n$ for some word $w_n\in X_{n+1}^*$.
\item \label{A6} For all $n\in\N$, we have $R_n\sub\Om$.
\item \label{A7} For all $n\in\N$, there is a mapping $X_n\to X_{n+1}^*:x\mt x^+$, and $\Om$ contains the relations
\[
x\lam_n = \lam_nx^+ \AND \rho_nx = x^+\rho_n \qquad\text{for all $x\in X_n$.}
\]
\item \label{A8} For all $n\in\N$, and for all $u\in X_{n+1}^*$, we have $w_nuw_n\sim w_nv^+w_n$ for some $v\in X_n^*$, where $w_n$ is the word from \ref{A5}, and $X_n^*\to X_{n+1}^*:u\mt u^+$ is (the extension of) the map from~\ref{A7}.
\een
The following is \cite[Theorem 2.16]{east2020presentations}, but see also \cite[Remark 2.12]{east2020presentations}:

\begin{theorem}\label{thm:Cpres}
With the above notation, and subject to assumptions \emph{\ref{A1}--\ref{A8}}, the category $\C$ has presentation $\Cpres\Ga\Om$ via $\phi$.  \epfres
\end{theorem}

The next result, also from \cite{east2020presentations}, shows how to convert the presentation from Theorem \ref{thm:Cpres} into a tensor presentation.  In what follows, $\C$ still denotes a tensor category satisfying the assumptions \ref{A1}--\ref{A8}, and we retain the meaning of $\Ga$, $\Om$, $\phi$ and so on.  We now also fix a digraph $\De$ on vertex set $\N$, a set of relations $\Xi\sub\De\fTC\times\De\fTC$, and a morphism $\Phi:\De\fTC\to\C$.  For $w\in\De\fTC$ we write $\ul w=w\Phi\in\C$, and we denote by $\approx$ the congruence $\Xi\tsharp$ on $\De\fTC$.  We also assume that the following all hold:
\begin{enumerate}[label=\textup{(A\arabic*)}]\addtocounter{enumi}{8}
\item \label{A9} For every relation $(u,v)\in\Xi$, we have $\ul u=\ul v$.
\item \label{A10} There is a mapping $\Ga\to\De\fTC:x\mt\wh x$ such that $\ulwh{x}=\ol x$ (i.e., $\wh x\Phi=x\phi$) for all $x\in\Ga$.  We extend this to a morphism $\Ga^*\to\De\fTC:w\mt\wh w$.  It quickly follows that $\ulwh{w}=\ol w$ for all $w\in\Ga^*$.
\item \label{A11} For a generator $x\in\De$, and for natural numbers $m,n\in\N$, we define the term
\[
x_{m,n} = \io_m\op x\op\io_n\in\De\fTC.
\]
We assume that for every such $x,m,n$, we have $x_{m,n}\approx\wh w$ for some $w\in\Ga^*$.
\item \label{A12} For every relation $(u,v)\in\Om$, we have $\wh u\approx\wh v$.
\een
The following is \cite[Theorem 2.17]{east2020presentations}:

\begin{theorem}\label{thm:TCpres}
With the above notation, and subject to assumptions \emph{\ref{A1}--\ref{A12}}, the category~$\C$ has tensor presentation $\TCpres\De\Xi$ via $\Phi$.  \epfres
\end{theorem}

\section{\boldmath The monoid $M\wr\I_n$}\label{sect:MIn}

Our goal in this section is to give monoid presentations for the wreath product
\[
M\wr\I_n = \bigset{(\ba,\al)}{\ba\in M_0^n,\ \al\in\I_n,\ \supp(\ba)=\dom(\al)}.
\]
We give two such presentations (Theorems \ref{thm:MIn} and \ref{thm:MIn2}), each extending a known presentation for~$\I_n$ (stated in Theorems \ref{thm:In} and \ref{th:popova}).  The first presentation for $M\wr\I_n$ is more symmetrical, and will be used in Section \ref{sec:category} when we study the category $M\wr\I$.  The second presentation utilises a smaller generating set, and will be obtained by re-writing the first. Examples \ref{eg:explicit} and \ref{eg:explicitsimp} show how the  presentations given in Theorem \ref{thm:MIn} and Theorem \ref{thm:MIn2} respectively look for small, explicit choices of $M$ and $n$.

\subsection[Presentations for $\I_n$]{\boldmath Presentations for $\I_n$}

As just noted, our presentations for $M\wr\I_n$ will utilise presentations for the symmetric inverse monoid $\I_n$, which we now state.  We begin by defining the alphabet
\begin{equation}\label{eq:XIn}
X_{\I_n} = \{s_1,\ldots,s_{n-1}\}\cup\{e_1,\ldots,e_n\},
\end{equation}
and the morphism
\[
\phi_{\I_n}:X_{\I_n}^*\to\I_n \qquad\text{given by}\qquad s_i\mt \ol s_i=(i,i+1) \AND e_i\mt\ol e_i = \trans{1&\cdots&i-1&i&i+1&\cdots&n\\1&\cdots&i-1&-&i+1&\cdots&n}.
\]
(Here $(i,i+1)$ is the transposition swapping $i$ and $i+1$.)
Diagrammatic representations of the $\ol s_i$ and $\ol e_i$ are given in Figure \ref{fig:minvmongens}.
We denote by $R_{\I_n}$ the following set of relations over $X_{\I_n}$, where in each relation the subscripts range over all meaningful values, subject to any stated constraints:
\begin{align}
\label{In1} s_i^2 &= \io , \\
\label{In2} s_is_j &= s_js_i &&\hspace{-2cm}\text{if $|i - j| > 1$,} \\
\label{In3} s_is_js_i &= s_js_is_j &&\hspace{-2cm}\text{if $|i - j| = 1$,} \\
\label{In4} e_i^2 &= e_i, \\
\label{In5} e_ie_j &= e_je_i, \\
\label{In6} s_ie_j &= e_js_i  &&\hspace{-2cm}\text{if $j\not=i,i+1$,} \\
\label{In7} s_ie_i &= e_{i+1}s_i, \\
\label{In8} e_ie_{i+1}s_i &= e_ie_{i+1}.
\end{align}
The next result follows quickly from \cite[Theorem 4.8]{Gilbert2006} and \cite[Proposition 31]{East2007}:

\begin{theorem}\label{thm:In}
For any $n\geq0$, the symmetric inverse monoid $\I_n$ has presentation ${\Mpres{X_{\I_n}}{R_{\I_n}}}$ via $\phi_{\I_n}$.  \epfres
\end{theorem}

It will also be convenient to state Popova's original presentation for $\I_n$ from \cite{popova1961defining}, which utilises a smaller generating set.  For this, we define the alphabet
\begin{equation}\label{eq:XIn'}
X_{\I_n}' = \{s_1,\ldots,s_{n-1},e\},
\end{equation}
and morphism
\[
\phi_{\I_n}':(X_{\I_n}')^*\to\I_n \qquad\text{given by}\qquad s_i\mt \ol s_i=(i,i+1) \AND e\mt\ol e = \trans{1&2&\cdots&n\\-&2&\cdots&n}.
\]
We also denote by $R_{\I_n}'$ the set of relations consisting of \eqref{In1}--\eqref{In3} together with
\begin{align}
\label{In9} e^2&=e, \\
\label{In10} es_1es_1=es_1e&=s_1es_1e,\\
\label{In11} es_i&=s_ie &&\hspace{-2cm}\text{if $i\geq 2$}.
\end{align}
The following is attributed in \cite{EL2004} to Popova \cite{popova1961defining}:

\begin{theorem}\label{th:popova}
For any $n\geq1$, the symmetric inverse monoid $\I_n$ has presentation ${\Mpres{X_{\I_n}'}{R_{\I_n}'}}$ via $\phi_{\I_n}'$.  \epfres
\end{theorem}

It should be noted that a symmetrical version of Popova's presentation is often stated, where the generator~$e$ maps instead to $\trans{1&\cdots&n-1&n\\1&\cdots&n-1&-}$.  One then has to modify the relations accordingly, so for example~\eqref{In10} becomes $es_{n-1}es_{n-1}=es_{n-1}e=s_{n-1}es_{n-1}e$.

\subsection[First presentation for $M\wr\I_n$]{\boldmath First presentation for $M\wr\I_n$}\label{subsect:MIn}

Our first presentation for $M\wr\I_n$ will extend the presentation for $\I_n$ from Theorem \ref{thm:In}, and will incorporate a presentation for $M$.  Thus, for the rest of this section we fix a presentation $\Mpres{X_M}{R_M}$ for~$M$, via $\phi_M:X_M^*\to M$.  For $w\in X_M^*$ we write $\ol w=w\phi_M$.

We begin by defining $n$ disjoint copies of $X_M$:
\[
X_M^{(i)} = \set{x^{(i)}}{x\in X_M} \qquad\text{for each $i\in\bn$,}
\]
and we set
\[
X_{M\wr\I_n} = X_{\I_n} \cup X_M^{(1)}\cup\cdots\cup X_M^{(n)},
\]
where $X_{\I_n}$ is as in \eqref{eq:XIn}.  

In order to define a morphism $X_{M\wr\I_n}^*\to M\wr\I_n$, we first introduce a piece of notation.  For $a\in M$ and $i\in\bn$, we define the tuple $a^{(i)} = (1,\ldots,1,a,1,\ldots,1)\in M^n$, where $a$ is in position~$i$.  We then define the morphism
\[
\phi_{M\wr\I_n}:X_{M\wr\I_n}^*\to M\wr\I_n : \begin{cases}
x \mt \ol x &\text{for $x\in X_{\I_n}$,}\\
x^{(i)} \mt \ol x^{(i)} &\text{for $x\in X_M$ and $i\in\bn$,}
\end{cases}
\]
where again we identify both $M^n$ and $\I_n$ as submonoids of $M\wr\I_n$, as explained at the end of Subsection \ref{subsect:wreath}.
Diagrammatic representations of the images of the letters from $X_{M\wr\I_n}$ are shown in Figure \ref{fig:minvmongens}.  The reader is reminded of the convention regarding suppressed upper labels.
As with words in $X_M^*$, we extend the over-line notation for words $w\in X_{M\wr\I_n}^*$ by writing $\ol w = w\phi_{M\wr\I_n}$.

\begin{figure}[ht]
\begin{center}
\begin{tikzpicture}[scale=.47]
\begin{scope}[shift={(0,0)}]	
\udotted14
\ddotted14
\udotted7{10}
\ddotted7{10}
\stline11
\stline44
\stline65
\stline56
\stline77
\stline{10}{10}
\node()at(1,2.6){\small$1$};
\node()at(5,2.6){\small$i$};
\node()at(10,2.6){\small$n$};
\uvs{1,4,5,6,7,10}
\lvs{1,4,5,6,7,10}
\end{scope}
\begin{scope}[shift={(13,0)}]	
\uvs{1,4,5,6,9}
\lvs{1,4,5,6,9}
\udotted14
\ddotted14
\udotted69
\ddotted69
\stline11
\stline44
\stline66
\stline{9}{9}
\node()at(1,2.6){\small$1$};
\node()at(5,2.6){\small$i$};
\node()at(9,2.6){\small$n$};
\end{scope}
\begin{scope}[shift={(25,0)}]	
\uvs{1,4,5,6,9}
\lvs{1,4,5,6,9}
\udotted14
\ddotted14
\udotted69
\ddotted69
\stline11
\stline44
\stline55
\stline66
\stline{9}{9}
\node()at(1,2.6){\small$1$};
\node()at(5,2.9){\small$i$};
\node()at(9,2.6){\small$n$};
\node[circle,draw=black, fill=white, inner sep = 0.06cm] () at (5,2){\footnotesize$\ol x$};
\end{scope}
\end{tikzpicture}
\caption{Generators of $M \wr \invmon$ (left to right): $\ol s_i$, $\ol e_i$ and $\ol x^{(i)}$.}
\label{fig:minvmongens}
\end{center}
\end{figure}

To state the relations for our presentation of $M\wr\I_n$ we need another piece of notation.  For any word ${w=x_1\cdots x_k\in X_M^*}$ and for any $i\in\bn$, we define the word
\[
w^{(i)} = x_1^{(i)}\cdots x_k^{(i)}\in(X_M^{(i)})^*.
\]
We now define $R_{M\wr\I_n}$ to be the set of relations consisting of $R_{\I_n}$ (i.e., \eqref{In1}--\eqref{In8}) together with the following, again with sub/superscripts $i$ and $j$ as well as elements $x$ and $y$ of $X_M$ ranging over all meaningful values, subject to any stated constraints:
\begin{align}
\label{MIn1} u^{(i)} &= v^{(i)} &&\hspace{-1cm}\text{for all $(u,v)\in R_M$ and $i\in\bn$,}\\
\label{MIn2} x^{(i)}y^{(j)} &= y^{(j)}x^{(i)} &&\hspace{-1cm}\text{if $i\not=j$,}\\
\label{MIn3} s_ix^{(j)} &= x^{(j)}s_i &&\hspace{-1cm}\text{if $j\not=i,i+1$,}\\
\label{MIn4} s_ix^{(i)} &= x^{(i+1)}s_i, \\
\label{MIn5} e_ix^{(j)} &= x^{(j)}e_i &&\hspace{-1cm}\text{if $i\not=j$,}\\
\label{MIn6} e_ix^{(i)} = e_i &= x^{(i)}e_i.
\end{align}

\begin{theorem}\label{thm:MIn}
For $n\geq0$, the wreath product $M\wr\I_n$ has presentation $\Mpres{X_{M\wr\I_n}}{R_{M\wr\I_n}}$ via $\phi_{M\wr\I_n}$.
\end{theorem}

\pf
To show that $\phi_{M\wr\I_n}$ is surjective, let $(\ba,\al)\in M\wr\I_n$.  Let $\bb=(b_1,\ldots,b_n)\in M^n$ be the tuple obtained from $\ba$ by replacing any $0$ entries by $1$, so that $(\ba,\al) = \bb\cdot\al$ (cf.~Remark \ref{rem:aabb}), and note that $\bb=b_1^{(1)}\cdots b_n^{(n)}$.  For each $i\in\bn$, we have $b_i = \ol w_i$ for some $w_i\in X_M^*$, and we also have $\al=\ol w$ for some $w\in X_{\I_n}^*$.  It then follows that
\[
(\ba,\al) = b_1^{(1)}\cdots b_n^{(n)}\cdot\al = \ol{w_1^{(1)}\cdots w_n^{(n)}\cdot w}.
\]
Next, it is easy to check that $\phi_{M\wr\I_n}$ preserves all the relations from $R_{M\wr\I_n}$; a couple of sample calculations are given in Figure~\ref{fig:MIn_rels} for \eqref{MIn4} and \eqref{MIn6}.  It follows from this that ${R_{M\wr\I_n}^\sharp \sub \ker(\phi_{M\wr\I_n})}$.
The rest of the proof is devoted to establishing the reverse inclusion.  For this, we write $\sim$ for the congruence $R_{M\wr\I_n}^\sharp$.  We begin with a technical lemma:

\begin{lemma}\label{lem:w1wnw'}
For any $w\in X_{M\wr\I_n}^*$, we have $w\sim w_1^{(1)}\cdots w_n^{(n)}\cdot w'$ for some $w'\in X_{\I_n}^*$, and some $w_1,\ldots,w_n\in X_M^*$ with $w_i=\io$ for any $i\in\bn\sm\dom(\ol w')$.
\end{lemma}

\pf
We first claim that
\[
w\sim w''w' \qquad\text{for some $w''\in(X_M^{(1)}\cup\cdots\cup X_M^{(n)})^*$ and $w'\in X_{\I_n}^*$.}
\]
Indeed, this follows from Lemma \ref{lem:XY}\ref{XY1}.  Note that to apply this lemma, we need to check that for any $x\in X_M^{(1)}\cup\cdots\cup X_M^{(n)}$ and $y\in X_{\I_n}$, we have $yx\sim uv$ for some $u\in(X_M^{(1)}\cup\cdots\cup X_M^{(n)})^*$ and $v\in X_{\I_n}^*$ with $\ell(u)\leq1$.  Now, $R_{M\wr\I_n}$ itself contains such a relation $(yx,uv)$ in every case except $x=x^{(i+1)}$ and $y=s_i$.  But for this we use \eqref{In1} and \eqref{MIn4} to calculate
\[
yx = s_ix^{(i+1)} \sim s_ix^{(i+1)}s_is_i \sim s_is_ix^{(i)}s_i \sim x^{(i)}s_i,
\]
and we can take $u=x^{(i)}$ and $v=s_i$.  This completes the proof of the claim, and we now fix~$w'$ and~$w''$ as above.

Next we use \eqref{MIn2} to `unshuffle' $w''$: i.e., to move all the letters in $w''$ from $X_M^{(1)}$ to the left, followed by all those from $X_M^{(2)}$, and so on.  Thus, we have $w''\sim w_1^{(1)}\cdots w_n^{(n)}$ for some ${w_1,\ldots,w_n\in X_M^*}$, from which it follows by our first claim that $w\sim w_1^{(1)}\cdots w_n^{(n)}\cdot w'$.

It remains to show that any $w_i^{(i)}$ with $i\not\in\dom(\ol w')$ can be deleted (i.e., replaced by $\io$).  To do so, fix some such $i$.  We then have $\ol w' = \ol e_i\ol w'$ (in $\I_n$), and so $w'\sim e_iw'$ by Theorem \ref{thm:In}.  Combining this with \eqref{MIn5} and \eqref{MIn6}, and writing $u=w_1^{(1)}\cdots w_{i-1}^{(i-1)}$ and $v=w_{i+1}^{(i+1)}\cdots w_n^{(n)}$ for convenience, we then have
\begin{align*}
w \sim u\cdot w_i^{(i)}\cdot v\cdot w' 
\sim u\cdot w_i^{(i)}\cdot v\cdot e_iw' 
\sim u\cdot w_i^{(i)}e_i\cdot v\cdot w' 
\sim u\cdot e_i\cdot v\cdot w' 
\sim u\cdot v\cdot e_iw' 
\sim u\cdot v\cdot w',
\end{align*}
completing the proof.
\epf

Returning now to the main proof, let $(u,v)\in\ker(\phi_{M\wr\I_n})$, so that $u,v\in X_{M\wr\I_n}^*$ and $\ol u=\ol v$; we must show that $u\sim v$.  By Lemma \ref{lem:w1wnw'}, we have
\[
u \sim u_1^{(1)}\cdots u_n^{(n)}\cdot u' \AND v \sim v_1^{(1)}\cdots v_n^{(n)}\cdot v',
\]
for appropriate $u_i,v_i\in X_M^*$ and $u',v'\in X_{\I_n}^*$.  Define $\ba,\bb\in M_0^n$ by
\[
a_i = \begin{cases}
\ol u_i &\text{if $i\in\dom(\ol u')$}\\
0 &\text{otherwise}
\end{cases}
\AND
b_i = \begin{cases}
\ol v_i &\text{if $i\in\dom(\ol v')$}\\
0 &\text{otherwise.}
\end{cases}
\]
Then, keeping Remark \ref{rem:aabb} in mind, we have
\[
((a_1,\ldots,a_n),\ol u') = (\ol u_1,\ldots,\ol u_n) \cdot \ol u' = \ol u = \ol v = (\ol v_1,\ldots,\ol v_n)\cdot \ol v' = ((b_1,\ldots,b_n),\ol v').
\]
In particular we have $\ol u'=\ol v'$, and it follows from Theorem \ref{thm:In} that $u'\sim v'$.  We also have $a_i=b_i$ for all $i$, so in particular when $i\in\dom(\ol u')=\dom(\ol v')$ we have $\ol u_i=\ol v_i$.  Since $M$ has presentation $\Mpres{X_M}{R_M}$, it follows that $(u_i,v_i)\in R_M^\sharp$ for any such $i$, and hence $u_i^{(i)} \sim v_i^{(i)}$ by \eqref{MIn1}.  When $i\not\in\dom(\ol u_i')=\dom(\ol v_i')$ we have $u_i=v_i=\io$, so also $u_i^{(i)}\sim v_i^{(i)}$ for such $i$.  It follows from all of this that
\[
u \sim u_1^{(1)}\cdots u_n^{(n)}\cdot u' \sim v_1^{(1)}\cdots v_n^{(n)}\cdot v' \sim v,
\]
and we have completed the proof of the theorem.
\epf 

\begin{figure}[ht]
\begin{center}
\begin{tikzpicture}[scale=.7]
\begin{scope}[shift={(0,0)}]	
\stline12
\stline21
\uvs{1,2}
\lvs{1,2}
\node () at (3.5,0) {$=$};
\end{scope}
\begin{scope}[shift={(0,-2)}]	
\stline11
\stline22
\uvs{1,2}
\lvs{1,2}
\node[circle,draw=black, fill=white, inner sep = 0.06cm] () at (1,2){\footnotesize$\ol x$};
\end{scope}
\begin{scope}[shift={(4,-2)}]	
\stline12
\stline21
\uvs{1,2}
\lvs{1,2}
\end{scope}
\begin{scope}[shift={(4,0)}]	
\stline11
\stline22
\uvs{1,2}
\lvs{1,2}
\node[circle,draw=black, fill=white, inner sep = 0.06cm] () at (2,2){\footnotesize$\ol x$};
\end{scope}
\begin{scope}[shift={(10,0)}]
\uvs{1}
\lvs{1}
\node () at (2.5,0) {$=$};
\end{scope}
\begin{scope}[shift={(10,-2)}]
\uvs{1}
\lvs{1}
\stline11
\node[circle,draw=black, fill=white, inner sep = 0.06cm] () at (1,2){\footnotesize$\ol x$};
\end{scope}
\begin{scope}[shift={(13,-1)}]
\uvs{1}
\lvs{1}
\node () at (2.5,1) {$=$};
\end{scope}
\begin{scope}[shift={(16,-2)}]
\uvs{1}
\lvs{1}
\end{scope}
\begin{scope}[shift={(16,0)}]
\uvs{1}
\lvs{1}
\stline11
\node[circle,draw=black, fill=white, inner sep = 0.06cm] () at (1,2){\footnotesize$\ol x$};
\end{scope}
\end{tikzpicture}
\caption{Left: relation \eqref{MIn4}, picturing only strings $i$ and $i+1$.  Right: relation \eqref{MIn6}, picturing only string $i$.}
\label{fig:MIn_rels}
\end{center}
\end{figure}
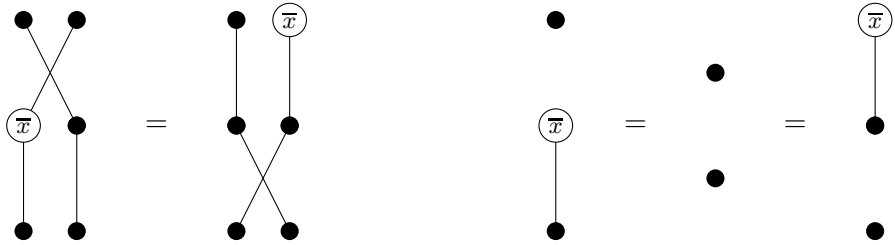

\begin{eg}\label{eg:explicit}
Before moving on it is worth considering an explicit, small example of the presentation $\Mpres{X_{M\wr\I_n}}{R_{M\wr\I_n}}$ from Theorem \ref{thm:MIn}. Specifically, we consider the case in which $n=2$ and $M=\pres{a,b}{ab=\io}$ is the bicyclic monoid.  Here we have
\[
X_{M\wr\I_2} = \{s_1,e_1,e_2,a^{(1)},b^{(1)},a^{(2)},b^{(2)}\}.
\]
Relations \eqref{In2}, \eqref{In3}, \eqref{In6} and \eqref{MIn3} are empty for $n=2$. The remaining relations from \eqref{In1}--\eqref{In8} have the form:
\bit\bmc2
\item[\eqref{In1}] $s_1^2=\io$, 
\item[\eqref{In4}] $e_1^2=e_1$, $e_2^2=e_2$, 
\item[\eqref{In5}] $e_1e_2=e_2e_1$, 
\item[\eqref{In7}] $s_1e_1=e_2s_1$, 
\item[\eqref{In8}] $e_1e_2s_1=e_1e_2$,
\item[] \phantom{a}
\emc\eit
while the remaining relations in \eqref{MIn1}--\eqref{MIn6} have the form:
\bit
\item[\eqref{MIn1}] $a^{(1)}b^{(1)} = a^{(2)}b^{(2)} = \io$,
\item[\eqref{MIn2}] 
$a^{(1)}a^{(2)}=a^{(2)}a^{(1)}$, 
$a^{(1)}b^{(2)}=b^{(2)}a^{(1)}$, 
$b^{(1)}a^{(2)}=a^{(2)}b^{(1)}$, 
$b^{(1)}b^{(2)}=b^{(2)}b^{(1)}$, 
\item[\eqref{MIn4}] $s_1a^{(1)}=a^{(2)}s_1$, $s_1b^{(1)}=b^{(2)}s_1$,
\item[\eqref{MIn5}] 
$e_1a^{(2)}=a^{(2)}e_1$, 
$e_1b^{(2)}=b^{(2)}e_1$, 
$e_2a^{(1)}=a^{(1)}e_2$, 
$e_2b^{(1)}=b^{(1)}e_2$, 
\item[\eqref{MIn6}] $e_1a^{(1)}=a^{(1)}e_1=e_1=e_1b^{(1)}=b^{(1)}e_1$, $e_2a^{(2)}=a^{(2)}e_2=e_2=e_2b^{(2)}=b^{(2)}e_2$.
\eit
\end{eg}

\subsection[Second presentation for $M\wr\I_n$]{\boldmath Second presentation for $M\wr\I_n$}\label{subsect:MIn2}

The presentation for $M\wr\I_n$ from Theorem \ref{thm:MIn} extended the presentation for $\I_n$ from Theorem~\ref{thm:In}.  We now give a second presentation, in terms of a smaller generating set, which will extended Popova's presentation for $\I_n$ from Theorem~\ref{th:popova}.

We continue to fix the presentation $\Mpres{X_M}{R_M}$ for $M$, via $\phi_M:X_M^*\to M$.  We begin by defining the alphabet
\[
X_{M\wr\I_n}' = X_{\I_n}' \cup X_M,
\]
where $X_{\I_n}'$ is as in \eqref{eq:XIn'}, and we define the morphism
\[
\phi_{M\wr\I_n}' : (X_{M\wr\I_n}')^* \to M\wr\I_n: \begin{cases}
x \mt \ol x &\text{for $x\in X_{\I_n}'$,}\\
x \mt (\ol{x}, 1, \dots, 1) &\text{for $x\in X_M$.}
\end{cases}
\]
Here we continue to identify $M^n$ and $\I_n$ with submonoids of $M\wr\I_n$.  Again we write $\ol w=w\phi_{M\wr\I_n}'$ for $w\in(X_{M\wr\I_n}')^*$.  Denote by $R_{M\wr\I_n}'$ the relations consisting of~$R_{\I_n}'$ (i.e., \eqref{In1}--\eqref{In3} and \eqref{In9}--\eqref{In11}), the relations in~$R_M$, as well as the following relations where $x, y \in X_M$:
\begin{align}
\label{MIn7} s_ix &= xs_i &&\hspace{-2cm}\text{if $i\geq2$,}\\
\label{MIn8} xs_1ys_1 &= s_1ys_1x ,\\
\label{MIn9} es_1xs_1 &= s_1xs_1e ,\\
\label{MIn10} ex=xe &=e.
\end{align}

\begin{theorem}\label{thm:MIn2}
For $n\geq1$, the wreath product $M\wr\I_n$ has presentation $\Mpres{X_{M\wr\I_n}'}{R_{M\wr\I_n}'}$ via $\phi_{M\wr\I_n}'$.
\end{theorem}

\pf
To simplify notation, throughout the proof we will write
\[
Y = X_{M\wr\I_n} \COMMA Z = X_{M\wr\I_n}' \COMMA S = R_{M\wr\I_n} \AND T= R_{M\wr\I_n}'.
\]
By Theorem \ref{thm:MIn} we know that $M\wr\I_n \cong \Mpres YS$, and we will prove the current theorem by showing that $\Mpres YS \cong \Mpres ZT$.  (This can be thought of as a variation of the Tietze transformation technique.)
For the rest of the proof we write ${\sim}$ for $S^\sharp$ and ${\approx}$ for $T^\sharp$.
For $u\in Y^*$ and $v\in Z^*$ we write $[u]$ and $\ldb v\rdb$ for the $\sim$-class of $u$ and the $\approx$-class of $v$, respectively.

We begin by defining two morphisms:
\[
\psi_1:Y^*\to Z^*:\begin{cases}
s_i\mt s_i,\\
e_i\mt s_{i-1}\cdots s_1es_1\cdots s_{i-1},\\
x^{(i)}\mt s_{i-1}\cdots s_1xs_1\cdots s_{i-1},
\end{cases}
 \AND
 \psi_2:Z^*\to Y^*:\begin{cases}
s_i\mt s_i,\\
e\mt e_1,\\
x\mt x^{(1)}.
\end{cases}
\]
The main substance of the proof involves demonstrating that the following all hold:
\ben
\item \label{MIn'1} For every relation $(u,v)\in S$, we have $u\psi_1\approx v\psi_1$.
\item \label{MIn'2} For every relation $(u,v)\in T$, we have $u\psi_2\sim v\psi_2$.
\item \label{MIn'3} For every generator $y\in Y$, we have $y\psi_1\psi_2\sim y$.
\item \label{MIn'4} For every generator $z\in Z$, we have $z\psi_2\psi_1\approx z$.
\een
Before we prove these, we first note that they imply the theorem.  Indeed, it follows from \ref{MIn'1} and~\ref{MIn'2} that we have well-defined morphisms
\[
\Psi_1:\Mpres YS \to \Mpres ZT \AND \Psi_2:\Mpres ZT\to\Mpres YS,
\]
given by $[u]\Psi_1=\ldb u\psi_1\rdb$ and $\ldb v\rdb\Psi_2=[v\psi_2]$, for $u\in Y^*$ and $v\in Z^*$.  By \ref{MIn'3} and \ref{MIn'4} we have $[y]\Psi_1\Psi_2=[y]$ and $\ldb z\rdb\Psi_2\Psi_1=\ldb z\rdb$ for all $y\in Y$ and $z\in Z$.  It follows from these that $\Psi_1\Psi_2$ and~$\Psi_2\Psi_1$ are identity maps, so that $\Psi_1$ is an isomorphism (with $\Psi_2=\Psi_1^{-1}$).

It now remains to prove items \ref{MIn'1}--\ref{MIn'4}, and we cover these roughly in order of difficulty.

\pfitem{\ref{MIn'4}}  This is trivial, as we actually have equality: $z\psi_2\psi_1= z$ for all $z\in Z$.

\pfitem{\ref{MIn'2}}  Fix some $(u,v)\in T$.  Because of Theorem \ref{thm:MIn}, we can show that $u\psi_2\sim v\psi_2$ by showing that $\phi_{M\wr\I_n}$ maps $u\psi_2$ and $v\psi_2$ to the same element of~$M\wr\I_n$.  This can be verified diagrammatically for each such relation $(u,v)$, similarly to Figure~\ref{fig:MIn_rels}.

\pfitem{\ref{MIn'3}}  By Theorem \ref{thm:MIn}, this again just amounts to a diagrammatic check that $\phi_{M\wr\I_n}$ maps $y\psi_1\psi_2$ and $y$ to the same element of $M\wr\I_n$.

\pfitem{\ref{MIn'1}}  We treat the relations from $S=R_{M\wr\I_n}$ (i.e., \eqref{In1}--\eqref{In8} and \eqref{MIn1}--\eqref{MIn6}) one at a time.  In what follows, we will write $\wh u=u\psi_1$ for all $u\in Y^*$.  We must show that $\wh u\approx\wh v$ for every relation $(u,v)\in S$.

\pfitem{\eqref{In1}--\eqref{In8}}  This follows immediately from the fact that $\I_n$ has presentations $\Mpres{X_{\I_n}}{R_{\I_n}}$ and $\Mpres{X_{\I_n}'}{R_{\I_n}'}$, by Theorems \ref{thm:In} and \ref{th:popova}, noting that $R_{\I_n}'\sub T$.

\pfitem{\eqref{MIn1}}  Fix some $(u,v)\in R_M$ and $i\in\bn$.  Writing $u=x_1\cdots x_k$, and using \eqref{In1}, we have
\begin{align*}
\wh{u^{(i)}} = \wh{x_1^{(i)}}\cdots\wh{x_k^{(i)}} 
&= (s_{i-1}\cdots s_1x_1s_1\cdots s_{i-1})\cdots(s_{i-1}\cdots s_1x_ks_1\cdots s_{i-1}) \\
& \approx s_{i-1}\cdots s_1 (x_1\cdots x_k) s_1\cdots s_{i-1} \\
&= s_{i-1}\cdots s_1 \cdot u \cdot s_1\cdots s_{i-1}.
\end{align*}
(If $u=\io$, then the previous conclusion still holds because of \eqref{In1}.)
A similar calculation gives $\wh{v^{(i)}} \approx s_{i-1}\cdots s_1\cdot v\cdot s_1\cdots s_{i-1}$.  Since $R_M\sub T$, it follows that 
\[
\wh{u^{(i)}} \approx s_{i-1}\cdots s_1 \cdot u \cdot s_1\cdots s_{i-1} \approx s_{i-1}\cdots s_1 \cdot v \cdot s_1\cdots s_{i-1} \approx \wh{v^{(i)}}.
\]

\pfitem{\eqref{MIn2}}  First consider some word $w\in\{s_1,\ldots,s_{n-1}\}^*$, and let $x\in X_M$ and $i\in\bn$.  We denote by $w^{-1}$ the reverse of the word $w$: i.e., if $w=s_{i_1}\cdots s_{i_k}$ then $w^{-1}=s_{i_k}\cdots s_{i_1}$.  Noting that $\ol w=w\phi_{\I_n}'$ is a permutation of $\bn$, we claim that
\begin{equation}\label{eq:xiw}
\wh{x^{(i)}}\cdot w \approx w \cdot \wh{x^{(i\ol w)}} \AND w^{-1}\cdot \wh{x^{(i)}}\cdot w \approx \wh{x^{(i\ol w)}}.
\end{equation}
By induction, it suffices to prove \eqref{eq:xiw} in the case that $w=s_k$ is a single letter.  For this, and given~\eqref{In1}, it suffices to show that
\[
s_k\cdot \wh{x^{(i)}}\cdot s_k \approx \begin{cases}
\wh{x^{(i)}} &\text{if $k\not=i,i-1$,}\\
\wh{x^{(i+1)}} &\text{if $k=i$,}\\
\wh{x^{(i-1)}} &\text{if $k=i-1$.}
\end{cases}
\]
When $k>i$ the result follows quickly from relations \eqref{In2} and \eqref{MIn7}, which say that $s_k$ commutes with each letter appearing in $\wh{x^{(i)}}$.  The $k=i$ case holds by equality, and the $k=i-1$ case by~\eqref{In1}.  For the final case, fix $k<i-1$.  We first note that 
\[
s_1\cdots s_{i-1}\cdot s_k\approx s_{k+1}\cdot s_1\cdots s_{i-1} \AND s_k\cdot s_{i-1}\cdots s_1\approx s_{i-1}\cdots s_1\cdot s_{k+1}.
\]
Indeed, these follow quickly from \eqref{In2} and \eqref{In3}, or by a diagrammatic check in~$\S_n$.  Using the relations \eqref{MIn7} and \eqref{In1} we then calculate
\begin{align*}
s_k\cdot \wh{x^{(i)}}\cdot s_k  = s_k\cdot s_{i-1}\cdots s_1xs_1\cdots s_{i-1}\cdot s_k &\approx s_{i-1}\cdots s_1\cdot s_{k+1}xs_{k+1}\cdot s_1\cdots s_{i-1}\\
&\approx s_{i-1}\cdots s_1\cdot s_{k+1}s_{k+1}x\cdot s_1\cdots s_{i-1}\\
&\approx s_{i-1}\cdots s_1\cdot x\cdot s_1\cdots s_{i-1} = \wh{x^{(i)}}.
\end{align*}
This completes the proof of \eqref{eq:xiw}.

Returning now to relation \eqref{MIn2}, fix some $x,y\in X_M$ and some $i,j\in\bn$ with $i\not=j$.  We must show that $\wh{x^{(i)}}\cdot \wh{y^{(j)}} \approx \wh{y^{(j)}}\cdot \wh{x^{(i)}}$.  Also fix some $w\in\{s_1,\ldots,s_{n-1}\}^*$ with $i=1\ol w$ and $j=2\ol w$.  Then by \eqref{eq:xiw} and \eqref{In1} we have
\[
\wh{x^{(i)}}\cdot \wh{y^{(j)}} = \wh{x^{(1\ol w)}}\cdot\wh{y^{(2\ol w)}} \approx w^{-1} \cdot \wh{x^{(1)}} \cdot w \cdot w^{-1} \cdot \wh{y^{(2)}} \cdot w \approx w^{-1} \cdot \wh{x^{(1)}} \cdot \wh{y^{(2)}} \cdot w,
\]
and similarly $\wh{y^{(j)}}\cdot \wh{x^{(i)}} \approx w^{-1} \cdot \wh{y^{(2)}} \cdot \wh{x^{(1)}} \cdot w$.  We then obtain $\wh{x^{(i)}}\cdot \wh{y^{(j)}} \approx \wh{y^{(j)}}\cdot \wh{x^{(i)}}$ from \eqref{MIn8}, which says precisely $\wh{x^{(1)}}\cdot \wh{y^{(2)}} \approx \wh{y^{(2)}}\cdot \wh{x^{(1)}}$.

\pfitem{\eqref{MIn3} and \eqref{MIn4}}  These are special cases of \eqref{eq:xiw}.

\pfitem{\eqref{MIn5}}  Similarly to \eqref{eq:xiw}, we can show that
\begin{equation}\label{eq:eiw}
\wh e_i\cdot w \approx w\cdot \wh e_{i\ol w} \AND \wh w^{-1}\cdot e_i\cdot w \approx \wh e_{i\ol w} \qquad\text{for all $i\in\bn$ and $w\in\{s_1,\ldots,s_{n-1}\}^*$.}
\end{equation}
(This also follows by a diagrammatic check, given that $\I_n$ has presentation $\Mpres{X_{\I_n}'}{R_{\I_n}'}$.)  But then for any $x\in X_M$ and $i,j\in\bn$ with $i\not=j$, we fix some $w\in\{s_1,\ldots,s_{n-1}\}^*$ such that $i=1\ol w$ and $j=2\ol w$, and we have
\[
\wh e_i \cdot \wh{x^{(j)}} \approx w^{-1} \cdot \wh e_1\cdot w \cdot w^{-1} \cdot \wh{x^{(2)}}\cdot w \approx w^{-1} \cdot \wh e_1\cdot \wh{x^{(2)}}\cdot w
\ANDSIM \wh{x^{(j)}} \cdot \wh e_i \approx w^{-1} \cdot \wh{x^{(2)}}\cdot \wh e_1 \cdot w.
\]
We then obtain $\wh e_i \cdot \wh{x^{(j)}} \approx \wh{x^{(j)}} \cdot \wh e_i$ from \eqref{MIn9}, which says that $\wh e_1\cdot \wh{x^{(2)}}\approx \wh{x^{(2)}}\cdot \wh e_1$.

\pfitem{\eqref{MIn6}}  This is similar to the previous case, but we use relation \eqref{MIn10} in place of \eqref{MIn9}.  
\epf

\begin{eg}\label{eg:explicitsimp}
As in Example \ref{eg:explicit}, we now consider the presentation $\Mpres{X_{M\wr\I_n}'}{R_{M\wr\I_n}'}$ from Theorem \ref{thm:MIn2} where $M =\pres{a,b}{ab=\io}$ and $n = 2$. Here we have
\[
X_{M\wr\I_2}' = \{s_1,e,a,b\}.
\]
This time, $R_{\I_2}'$ consists of the relations $s_1^2=\io$, $e^2=e$ and $es_1es_1=es_1e=s_1es_1e$, $R_M$ consists of the relation $ab=\io$, while relations \eqref{MIn7}--\eqref{MIn10} have the form:
\bit
\item[\eqref{MIn7}] empty for $n=2$,
\item[\eqref{MIn8}] $as_1as_1=s_1as_1a$, $as_1bs_1=s_1bs_1a$, $bs_1as_1=s_1as_1b$, $bs_1bs_1=s_1bs_1b$, 
\item[\eqref{MIn9}] $es_1as_1=s_1as_1e$, $es_1bs_1=s_1bs_1e$, 
\item[\eqref{MIn10}] $ea=ae=e=eb=be$.
\eit
\end{eg}

\section{\boldmath The category $M\wr\I$}\label{sec:category}

In this section we turn our attention to the category 
\[
M\wr\I = \bigset{(\ba,\al)\in \M\times\I}{\bd(\ba)=\bd(\al),\ \supp(\ba)=\dom(\al)}.
\]
Our main results here are Theorems \ref{thm:MI} and \ref{thm:MI2}, which give category and tensor category presentations for $M\wr\I$.  These theorems are proved by combining Theorems \ref{thm:Cpres} and \ref{thm:TCpres} (from~\cite{east2020presentations}) with Theorem \ref{thm:MIn}, which gives presentations for the endomorphism monoids of $M\wr\I$.

\subsection[Category presentation for $M\wr\I$]{\boldmath Category presentation for $M\wr\I$}\label{subsect:MI}

We begin by giving a category presentation for $M\wr\I$, using Theorem \ref{thm:Cpres}.  To apply this theorem, we require presentations for the endomorphism monoids $M\wr\I_n$ ($n\in\N$), and the most convenient such presentations are those from Theorem \ref{thm:MIn}.  However, to ensure that the generating sets are pairwise disjoint, the generators $s_i, e_i$ and $x^{(i)}$ from $X_{M\wr\I_n}$ for any $n \in \N$ will be renamed $s_{i;n}$, $e_{i;n}$ and $x^{(i;n)}$ respectively.
We again fix the presentation $\Mpres{X_M}{R_M}$ for $M$ via $\phi_M:X_M^*\to M$.  

With the above notational comments in mind, we begin by defining the digraph~$\Ga_{M\wr\I}$ over vertex set $\N$, and with edges
\bit
\item $s_{i;n}$ for each $n\in\N$ and $1\leq i<n$, with $\bd(s_{i;n})=\br(s_{i;n})=n$,
\item $e_{i;n}$ for each $n\in\N$ and $1\leq i\leq n$, with $\bd(e_{i;n})=\br(e_{i;n})=n$,
\item $x^{(i;n)}$ for each $x\in X_M$, $n\in\N$ and $1\leq i\leq n$, with $\bd(x^{(i;n)})=\br(x^{(i;n)})=n$,
\item $\lam_n,\rho_n$ for each $n\in\N$, with $\bd(\lam_n)=\br(\rho_n)=n$ and $\br(\lam_n)=\bd(\rho_n)=n+1$.
\eit
We define the morphism 
\[
\phi_{M\wr\I}:\Ga_{M\wr\I}^*\to M\wr\I:x\mt\ol x
\]
diagrammatically in Figure \ref{fig:phi}, which also pictures the identities $\ol\io_n$ ($n\in\N$).  (Here as usual $\io_n\in\Ga_{M\wr\I}^*$ denotes the empty path $n\to n$.)

\begin{figure}[ht]
\begin{center}
\begin{tikzpicture}[scale=.5]
\begin{scope}[shift={(16,0)}]	
\uvs{1,6}
\lvs{1,6,7}
\udotted16
\ddotted16
\stline11
\stline66
\draw(0.5,1)node[left]{$\ol\lam_n=$};
\node()at(1,2.5){\tiny$1$};
\node()at(6,2.5){\tiny$n$};
\end{scope}
\begin{scope}[shift={(16,-5)}]	
\uvs{1,6,7}
\lvs{1,6}
\udotted16
\ddotted16
\stline11
\stline66
\draw(0.5,1)node[left]{$\ol\rho_n=$};
\node()at(1,2.5){\tiny$1$};
\node()at(6,2.5){\tiny$n$};
\end{scope}
\begin{scope}[shift={(16,-10)}]	
\uvs{1,6}
\lvs{1,6}
\udotted16
\ddotted16
\stline11
\stline66
\draw(0.5,1)node[left]{$\ol\io_n=$};
\node()at(1,2.5){\tiny$1$};
\node()at(6,2.5){\tiny$n$};
\end{scope}
\begin{scope}[shift={(0,0)}]	
\uvs{1,4,5,6,7,10}
\lvs{1,4,5,6,7,10}
\udotted14
\ddotted14
\udotted7{10}
\ddotted7{10}
\stline11
\stline44
\stline56
\stline65
\stline77
\stline{10}{10}
\draw(0.5,1)node[left]{$\ol s_{i;n}=$};
\node()at(1,2.5){\tiny$1$};
\node()at(5,2.5){\tiny$i$};
\node()at(10,2.5){\tiny$n$};
\end{scope}
\begin{scope}[shift={(0,-5)}]	
\uvs{1,4,5,6,10}
\lvs{1,4,5,6,10}
\udotted14
\ddotted14
\udotted6{10}
\ddotted6{10}
\stline11
\stline44
\stline66
\stline{10}{10}
\draw(0.5,1)node[left]{$\ol e_{i;n}=$};
\node()at(1,2.5){\tiny$1$};
\node()at(5,2.5){\tiny$i$};
\node()at(10,2.5){\tiny$n$};
\end{scope}
\begin{scope}[shift={(0,-10)}]	
\uvs{1,4,5,6,10}
\lvs{1,4,5,6,10}
\udotted14
\ddotted14
\udotted6{10}
\ddotted6{10}
\stline11
\stline44
\stline55
\stline66
\stline{10}{10}
\draw(0.5,1)node[left]{$\ol x^{(i;n)}=$};
\node()at(1,2.5){\tiny$1$};
\node()at(5,2.7){\tiny$i$};
\node()at(10,2.5){\tiny$n$};
\node[circle,draw=black, fill=white, inner sep = 0.06cm] () at (5,2){\footnotesize$\ol x$};
\end{scope}
\end{tikzpicture}
\caption{Generators of $M\wr\I$, as well as the identity $\ol\io_n$.}
\label{fig:phi}
\end{center}
\end{figure}

For a word $u=x_1\cdots x_k\in X_M^*$, and for $n\in\N$ and $1\leq i\leq n$, we define
\[
u^{(i;n)} = x_1^{(i;n)}\cdots x_k^{(i;n)}\in\Ga_{M\wr\I}^*.
\]

\begin{remark}\label{rem:ioin}
If $u$ is empty, we interpret $u^{(i;n)}=\io_n$.
\end{remark}

Now let $\Om_{M\wr\I}$ be the following set of relations over $\Ga_{M\wr\I}$, where as usual the sub/superscripts and elements of~$X_M$ range over all meaningful values, subject to any stated constraints:
\begin{align}
\label{MI1} & s_{i;n}^2 = \io_n \COMMA e_{i;n}^2=e_{i;n} \COMMA e_{i;n}e_{j;n}=e_{j;n}e_{i;n} \COMMA s_{i;n}e_{j;n}=e_{j;n}s_{i;n} \text{ if $j\not=i,i+1$,} \\
\label{MI2} & s_{i;n}e_{i;n}=e_{i+1;n}s_{i;n} \COMMA e_{i;n}e_{i+1;n}s_{i;n}=e_{i;n}e_{i+1;n},  \\
\label{MI3} & s_{i;n}s_{j;n} = s_{j;n}s_{i;n}  \text{ if $|i-j|>1$}\COMMA s_{i;n}s_{j;n}s_{i;n} = s_{j;n}s_{i;n}s_{j;n} \text{ if $|i-j|=1$,} \\
\label{MI4} &u^{(i;n)} = v^{(i;n)} \text{ for $(u,v)\in R_M$} \COMMA x^{(i;n)}y^{(j;n)} = y^{(j;n)}x^{(i;n)} \text{ if $i\not=j$} ,\\
\label{MI5} &s_{i;n}x^{(j;n)} = x^{(j;n)}s_{i;n} \text{ if $j\not=i,i+1$} \COMMA s_{i;n}x^{(i;n)} = x^{(i+1;n)}s_{i;n},\\
\label{MI6} &e_{i;n}x^{(j;n)} = x^{(j;n)}e_{i;n} \text{ if $i\not=j$}\COMMA e_{i;n}x^{(i;n)}=e_{i;n}=x^{(i;n)}e_{i;n},\\
\label{MI7} &\lam_n\rho_n=\io_n \COMMA \rho_n\lam_n=e_{n+1;n+1} ,\\
\label{MI8} &s_{i;n}\lam_n = \lam_ns_{i;n+1} \COMMA e_{i;n}\lam_n = \lam_ne_{i;n+1} \COMMA x^{(i;n)}\lam_n = \lam_nx^{(i;n+1)} ,\\
\label{MI9} &\rho_ns_{i;n}=s_{i;n+1}\rho_n \COMMA \rho_ne_{i;n}=e_{i;n+1}\rho_n \COMMA \rho_nx^{(i;n)}=x^{(i;n+1)} \rho_n.
\end{align}

\begin{theorem}\label{thm:MI}
The category $M\wr\I$ has presentation $\Cpres{\Ga_{M\wr\I}}{\Om_{M\wr\I}}$ via $\phi_{M\wr\I}$.
\end{theorem}

\pf
As we have already said, we will be applying Theorem \ref{thm:Cpres}.  To do so, we must check that each of assumptions \ref{A1}--\ref{A8} hold, as listed in Subsection \ref{subsect:pres}.  In the following we write $\ol w=w\phi_{M\wr\I}$ for $w\in\Ga_{M\wr\I}^*$, and we also write ${\sim}$ for $\Om_{M\wr\I}^\sharp$ -- the congruence on $\Ga_{M\wr\I}^*$ generated by the relations~$\Om_{M\wr\I}$.

\pfitem{\ref{A1}, \ref{A5} and \ref{A6}}   These are clear.

\pfitem{\ref{A2} and \ref{A4}}  These are easily checked diagrammatically (see Figure \ref{fig:phi} for \ref{A2}).

\pfitem{\ref{A3}}  We require a presentation $\Mpres{X_n}{R_n}$ for each endomorphism monoid $M\wr\I_n$ ($n\in\N$).  For this we take the presentation from Theorem \ref{thm:MIn}, modified to ensure that the alphabets $X_n$ are pairwise disjoint, as explained at the beginning of this subsection.  Explicitly, we have
\[
X_n = \{s_{1;n},\ldots,s_{n-1;n}\} \cup \{e_{1;n},\ldots,e_{n;n}\} \cup \set{x^{(i;n)}}{x\in X_M,\ i\in\bn},
\]
and we take $R_n$ to be the suitably modified versions of relations $R_{M\wr\I_n}$.  The union of these modified relations are precisely \eqref{MI1}--\eqref{MI6}.

\pfitem{\ref{A7}}  Here we have $s_{i;n}^+=s_{i;n+1}$, and similarly for $e_{i;n}$ and $x^{(i;n)}$; see \eqref{MI8} and \eqref{MI9}.

\pfitem{\ref{A8}}  Fix some $n\in\N$ and some $u\in X_{n+1}^*$.  We must show that
\begin{equation}\label{eq:ewe}
e_{n+1;n+1}ue_{n+1;n+1} \sim e_{n+1;n+1}v^+e_{n+1;n+1} \qquad\text{for some $v\in X_n^*$.}
\end{equation}
Here $X_n^*\to X_{n+1}^*:v\mt v^+$ is the morphism extending the map $X_n\to X_{n+1}^*$ from \ref{A7}.  For simplicity, write $e=e_{n+1;n+1}$, and consider the element $(\ba,\al) = \ol{eue}\in M\wr\I_{n+1}$.  Because of the terms $\ol e=\ol e_{n+1;n+1}$, we see that $\dom(\al)$, $\im(\al)$ and $\supp(\ba)$ are all contained in $\bn$.  Now consider $(\bb,\be)\in M\wr\I_n$, where~$\be$ is the restriction of $\al$ to $\bn$, and where $\bb$ consists of the first~$n$ entries of~$\ba$.  Since $(\bb,\be)\in M\wr\I_n$, we have $(\bb,\be) = \ol v$ for some word $v\in X_n^*$.  By construction, we have $(\ba,\al)=\ol{ev^+e}$.  So $\ol{eue}=\ol{ev^+e}$, and since $M\wr\I_{n+1}$ has presentation $\Mpres{X_{n+1}}{R_{n+1}}$, it follows that $eue \sim ev^+e$.  This completes the proof of \eqref{eq:ewe}, and hence of the theorem.
\epf

\subsection[Tensor presentation for $M\wr\I$]{\boldmath Tensor presentation for $M\wr\I$}\label{subsect:MI2}

We now use Theorem \ref{thm:TCpres} to convert the category presentation for $M\wr\I$ from Theorem \ref{thm:MI} into a tensor presentation.  

We continue to fix the presentation $\Mpres{X_M}{R_M}$ for $M$ via $\phi_M:X_M^*\to M$.  We begin by defining the digraph $\De_{M\wr\I}$ on vertex set $\N$ with edges $\{X,U,\U\}\cup X_M$, with sources and targets given by
\[
\bd(X)=\br(X)=2 \COMMa \bd(U)=\br(\U)=\bd(x)=\br(x)=1 \text{ for all $x\in X_M$} \ANd \br(U)=\bd(\U)=0.
\]
(Note here that `$X$' is the name of a single edge of $\De_{M\wr\I}$.)  We define the morphism
\[
\Phi_{M\wr\I}:\De_{M\wr\I}\fTC\to M\wr\I: \begin{cases}
X \mt \ul X = \ol s_{1;2} ,\\
U\mt \ul U = \ol\rho_0, \\
\U\mt \ul \U = \ol\lam_0, \\
x\mt\ul x = \ol x^{(1;1)},
\end{cases}
\]
using the notation of the previous subsection.  The images of the generators under $\Phi_{M\wr\I}$ are shown in Figure~\ref{fig:Phi}.  In what follows, it is convenient to write $I=\io_1$ for the empty path $1\to1$; note that $\ul I$ is also shown in Figure \ref{fig:Phi}.  For a term $w\in\De_{M\wr\I}\fTC$ we write $\ul w=w\Phi_{M\wr\I}\in M\wr\I$.

\begin{figure}[ht]
\begin{center}
\begin{tikzpicture}[scale=.5]
\begin{scope}[shift={(0,0)}]	
\uvs{1,2}
\lvs{1,2}
\stline12
\stline21
\draw(0.5,1)node[left]{$\ul X=$};
\end{scope}
\begin{scope}[shift={(8,0)}]	
\uvs{1}
\draw(0.5,1)node[left]{$\ul U=$};
\end{scope}
\begin{scope}[shift={(15,0)}]	
\lvs{1}
\draw(0.5,1)node[left]{$\ul \U=$};
\end{scope}
\begin{scope}[shift={(22,0)}]	
\lvs{1}
\uvs{1}
\draw(0.5,1)node[left]{$\ul x=$};
\stline11
\node[circle,draw=black, fill=white, inner sep = 0.06cm] () at (1,2){\footnotesize$\ol x$};
\end{scope}
\begin{scope}[shift={(29,0)}]	
\lvs{1}
\uvs{1}
\stline11
\draw(0.5,1)node[left]{$\ul I=$};
\end{scope}
\end{tikzpicture}
\caption[blah]{Generators of $M\wr\I$, as well as the identity $\ul I=\ol\io_1$.}
\label{fig:Phi}
\end{center}
\end{figure}

We now let $\Xi_{M\wr\I}$ be the set of relations over $\De_{M\wr\I}$ consisting of $R_M$ (see Remark \ref{rem:RM}) together with the following, where $x\in X_M$ in the third line:
\begin{align}
\label{MI10} &X\circ X=\io_2 \COMMA (X\op I)\circ(I\op X)\circ(X\op I) = (I\op X)\circ(X\op I)\circ(I\op X),\\
\label{MI11} &\U\circ U = \io_0  \COMMA X\circ(U\op I) = I\op U \COMMA (\U\op I)\circ X = I\op \U  ,\\
\label{MI12} &X\circ(x\op I) = (I\op x)\circ X \COMMA x\circ U=U \COMMA \U\circ x=\U.
\end{align}

\begin{remark}\label{rem:RM}
We have stipulated that $R_M\sub\Xi_{M\wr\I}$, and a clarification is in order here.  Any word $u=x_1\cdots x_k\in X_M^*$ involved in a relation from $R_M$ is identified with the term $u\equiv x_1\circ\cdots\circ x_k$ from~$\De_{M\wr\I}\fTC$.  If $u$ is empty, then this term is interpreted as the empty path $\io_1=I$.
\end{remark}

\begin{theorem}\label{thm:MI2}
The category $M\wr\I$ has tensor presentation $\TCpres{\De_{M\wr\I}}{\Xi_{M\wr\I}}$ via $\Phi_{M\wr\I}$.
\end{theorem}

\pf
The proof will be via an application of Theorem \ref{thm:TCpres}.  We keep the notation of Subsection~\ref{subsect:MI}, in particular $\Ga_{M\wr\I}$, $\Om_{M\wr\I}$, $\phi_{M\wr\I}$ and so on.  Throughout the proof we write ${\approx}$ for the tensor congruence $(\Xi_{M\wr\I})\tsharp$ on~$\De_{M\wr\I}\fTC$ generated by the relations $\Xi_{M\wr\I}$.  In the proof of Theorem~\ref{thm:MI} we checked that assumptions \ref{A1}--\ref{A8} hold, so it remains to check \ref{A9}--\ref{A12}.

\pfitem{\ref{A9}}  It is easy to check that $\ul u=\ul v$ in $M\wr\I$ for each relation $(u,v)$ from $\Xi_{M\wr\I}$.  For example, the left-hand diagram of Figure \ref{fig:MIn_rels} already does this for the first relation in \eqref{MI12}.  

\pfitem{\ref{A10}}  We define the mapping $\Ga_{M\wr\I}\to\De_{M\wr\I}\fTC$ as follows, for all appropriate values of $i,n\in\N$: 
\bit\bmc2
\item $\wh s_{i;n} = \io_{i-1}\op X\op \io_{n-i-1}$, 
\item $\wh e_{i;n} = \io_{i-1}\op U\op\U\op \io_{n-i}$, 
\item $\wh x^{(i;n)} = \io_{i-1}\op x\op \io_{n-i}$ for $x\in X_M$, 
\item $\wh \lam_n = \io_n\op\U$, 
\item $\wh \rho_n = \io_n\op U$. 
\item[] ~
\emc\eit
It is routine to check diagrammatically that $\ulwh{x}=\ol x$ (i.e., $\wh x\Phi_{M\wr\I}=x\phi_{M\wr\I}$) for all $x\in\Ga_{M\wr\I}$.  We extend this mapping to a morphism $\Ga_{M\wr\I}^*\to\De_{M\wr\I}\fTC:w\mt\wh w$, and note that $\wh\io_n=\io_n$ for all $n$.

\pfitem{\ref{A11}}  We must show that each term $x_{m,n}=\io_m\op x\op\io_n\in\De_{M\wr\I}\fTC$ ($x\in\De_{M\wr\I}$, $m,n\in\N$) is $\approx$-equivalent to some term of the form $\wh w$, with $w\in\Ga_{M\wr\I}^*$.  This is clear for $x=X$ or for $x\in X_M$, as 
\[
X_{m,n} = \io_m\op X\op\io_n = \wh s_{m+1;m+n+2} \AND x_{m,n} = \io_m\op x\op\io_n = \wh x_{m+1;m+n+1}.
\]
For the case $x=U$ ($x=\U$ is symmetrical) we have
\[
U_{m,n} \approx \wh\si_{m+1;m+n+1} \circ\cdots\circ \wh\si_{m+n;m+n+1} \circ\wh\rho_{m+n},
\]
which follows from a diagrammatic check and an application of Theorem \ref{thm:I}.  (Alternatively it can be proved directly in exactly the same way as \cite[Lemma 3.7]{east2020presentations}; the proof from \cite{east2020presentations} concerned a different category containing $\I$, but used only relations from $\Xi_{M\wr\I}$.)

\pfitem{\ref{A12}}  Consider a relation $(u,v)\in\Om_{M\wr\I}$ (i.e., one of \eqref{MI1}--\eqref{MI9}).  We must show that $\wh u\approx\wh v$.  In any case that the terms $\wh u$ and $\wh v$ involve only the generators $X$, $U$ and $\U$, this follows immediately from Theorem \ref{thm:I} and a simple diagrammatic check that $\wh u\Phi_{M\wr\I} = \wh v\Phi_{M\wr\I}$.  This takes care of all relations from $\Om_{M\wr\I}$ other than \eqref{MI4}--\eqref{MI6}, and the third parts of each of \eqref{MI8} and \eqref{MI9}.  We now consider these remaining relations.

\pfitem{\eqref{MI4}}  For the first part of this relation, let $(u,v)\in R_M$ and $1\leq i\leq n$; we must show that $\wh u^{(i;n)}\approx\wh v^{(i;n)}$.  Writing $u=x_1\cdots x_k$, we first use the tensor category axioms to calculate
\begin{align*}
\wh u^{(i;n)} = \wh x_1^{(i;n)}\circ\cdots\circ\wh x_k^{(i;n)} &= (\io_{i-1}\op x_1\op \io_{n-i}) \circ \cdots \circ (\io_{i-1}\op x_k\op \io_{n-i}) \\
&= (\io_{i-1}\circ\cdots\circ\io_{i-1}) \op (x_1\circ\cdots\circ x_k) \op (\io_{n-i}\circ\cdots\circ\io_{n-i}) 
\\&= \io_{i-1} \op u \op \io_{n-i}.
\end{align*}
(Keeping Remarks \ref{rem:ioin} and \ref{rem:RM} in mind, we note that if $u=\io_1$ is empty, then we still have $\wh u^{(i;n)} = \wh\io_n=\io_n=\io_{i-1}\op\io_1\op\io_{n-i}=\io_{i-1}\op u\op\io_{n-i}$.)
A similar calculation gives $\wh v^{(i;n)} = \io_{i-1} \op v \op \io_{n-i}$, so $\wh u^{(i;n)}\approx\wh v^{(i;n)}$ now follows from the fact that~$R_M\sub\Xi_{M\wr\I}$.

For the second part of \eqref{MI4}, let $x,y\in X_M$ and $i,j\in\bn$ with $i\not=j$; by symmetry we can assume that $i<j$.  We must show that $\wh x^{(i;n)}\circ\wh y^{(j;n)} \approx \wh y^{(j;n)}\circ\wh x^{(i;n)}$, and in fact we have equality using only the tensor category axioms.  Indeed, we have
\begin{align*}
\wh x^{(i;n)}\circ\wh y^{(j;n)} &= (\io_{i-1} \op x \op \io_{j-1-i} \op I\op \io_{n-j}) \circ (\io_{i-1} \op I \op \io_{j-1-i} \op y\op \io_{n-j}) \\
&= (\io_{i-1}\circ\io_{i-1}) \op (x\circ I) \op (\io_{j-1-i}\circ\io_{j-1-i}) \op (I\circ y) \op (\io_{n-j}\circ\io_{n-j}) \\
&= \io_{i-1} \op x \op \io_{j-1-i} \op y \op \io_{n-j},
\end{align*}
and similarly $\wh y^{(j;n)}\circ\wh x^{(i;n)} = \io_{i-1} \op x \op \io_{j-1-i} \op y \op \io_{n-j}$.

\pfitem{\eqref{MI5}}  The first part is treated in the same fashion as the second part of \eqref{MI4}, where we again need only the tensor category axioms.  For the second part of \eqref{MI5}, it follows quickly from the axioms that
\[
\wh s_{i;n}\circ\wh x^{(i;n)} = \io_{i-1} \op (X\circ(x\op I)) \op \io_{n-i-1} \AND \wh x^{(i+1;n)}\circ\wh s_{i;n} = \io_{i-1} \op ((I\op x)\circ X) \op \io_{n-i-1},
\]
so $\wh s_{i;n}\circ\wh x^{(i;n)} \approx \wh x^{(i+1;n)}\circ\wh s_{i;n}$ follows from the first relation in \eqref{MI12}.

\pfitem{\eqref{MI6}}  The first part is again similar to the second part of \eqref{MI4}.  For the second part of \eqref{MI6}, we have
\[
\wh e_{i;n}\circ \wh x^{(i;n)} = \io_{i-1}\op((U\op\U)\circ x)\op\io_{n-i} \AND \wh x^{(i;n)}\circ\wh e_{i;n} = \io_{i-1}\op(x\circ(U\op\U))\op\io_{n-i},
\]
so here we need to show that $(U\op\U)\circ x \approx U\op\U \approx x\circ(U\op\U)$.  But this follows quickly from the second and third parts of \eqref{MI12}, and the fact that $U\op\U = U\circ\U$ (more generally, once can easily check that $a\op b=a\circ b$ in any tensor category whenever $\br(a)=\bd(b)=0$).  Indeed,
\[
(U\op\U)\circ x = (U\circ\U)\circ x = U\circ(\U\circ x) \approx U\circ\U = U\op\U,
\]
and similarly $x\circ(U\op\U)\approx U\op\U$.

\pfitem{\eqref{MI8} and \eqref{MI9}}  As noted earlier, we need only consider the third relation in each, which amounts to showing that
\[
\wh x^{(i;n)}\circ\wh\lam_n \approx \wh\lam_n\circ\wh x^{(i;n+1)} \AND \wh\rho_n\circ\wh x^{(i;n)} \approx \wh x^{(i;n+1)}\circ\wh\rho_n.
\]
By symmetry it suffices to just deal with the first, and here again we have equality, as
\begin{align*}
\wh x^{(i;n)} \circ\lam_n = (\wh x^{(i;n)} \op \io_0) \circ (\io_n\op\U) = (\wh x^{(i;n)}\circ\io_n) \op (\io_0\circ\U) = \wh x^{(i;n)}\op\U,
\intertext{and}
\wh\lam_n\circ\wh x^{(i;n+1)} = (\io_n\op\U) \circ (\wh x^{(i;n)}\op I) = (\io_n\circ\wh x^{(i;n)}) \op (\U\circ I) = \wh x^{(i;n)} \op \U.
\end{align*}
This completes the proof.
\epf

\section{\boldmath The semigroup $M\wr\Sing(\I_n)$}\label{sec:wreathprodsing}

In this final section we give a semigroup presentation for the singular wreath product 
\[
M \wr \Sing(\I_n) = \bigset{(\ba,\al)\in M_0^n\times \Sing(\I_n)}{\supp(\ba)=\dom(\al)}.
\]
(Recall that $\Sing(\I_n)=\I_n\sm\S_n$ is the singular part of $\I_n$.)  Our presentation (see Theorem~\ref{thm:MSIn}) will be built out of the presentation for $\Sing(\I_n)$ from \cite{East2015}, stated below in Theorem \ref{thm:SI}, and a new presentation for a certain (singular) subsemigroup of $M_0^n$, which we will give in Theorem~\ref{thm:M}.  Since $M\wr\Sing(\I_0)=\es$ and $M\wr\Sing(\I_1)=\{\es\}$, we assume that $n\geq2$ from here on.

\subsection[Presentation for $\Sing(\I_n)$]{\boldmath Presentation for $\Sing(\I_n)$}

To state the presentation for $\Sing(\I_n)$ from \cite{East2015}, we begin by defining the alphabet
\begin{equation}\label{eq:XS}
\XS = \set{f_{i,j}}{i,j\in\bn,\ i\not=j},
\end{equation}
and the morphism
\begin{equation}\label{eq:pS}
\pS:\XS^+\to\Sing(\I_n):f_{i,j}\mt\ol f_{i,j},
\end{equation}
where for distinct $i,j\in\bn$, the partial map $\ol f_{i,j}\in\Sing(\I_n)$ is defined by
\[
k \ol f_{i,j} = \begin{cases}
k &\text{if $k\not=i,j$,}\\
i &\text{if $k=j$,}\\
- &\text{if $k=i$.}
\end{cases}
\]
The maps $\ol f_{i,j}$ are pictured in Figure \ref{fig:f}.  Note that $\dom(\ol f_{i,j}) = \{i\}\com$ and $\im(\ol f_{i,j}) = \{j\}\com$.  Here and throughout, for $A\sub\bn$ we write $A\com=\bn\sm A$ for the complement of~$A$ in $\bn$.

\begin{figure}[ht]
\begin{center}
\begin{tikzpicture}[scale=.5]
\begin{scope}[shift={(0,0)}]	
\udotted14
\ddotted14
\udotted69
\ddotted69
\udotted{11}{14}
\ddotted{11}{14}
\stline11
\stline44
\stline66
\stline99
\stline{10}5
\stline{11}{11}
\stline{14}{14}
\node()at(1,2.6){\small$1$};
\node()at(5,2.6){\small$i$};
\node()at(10,2.6){\small$j$};
\node()at(14,2.6){\small$n$};
\uvs{1,4,5,6,9,10,11,14}
\lvs{1,4,5,6,9,10,11,14}
\end{scope}
\begin{scope}[shift={(18,0)}]	
\udotted14
\ddotted14
\udotted69
\ddotted69
\udotted{11}{14}
\ddotted{11}{14}
\stline11
\stline44
\stline66
\stline99
\stline5{10}
\stline{11}{11}
\stline{14}{14}
\node()at(1,2.6){\small$1$};
\node()at(5,2.6){\small$j$};
\node()at(10,2.6){\small$i$};
\node()at(14,2.6){\small$n$};
\uvs{1,4,5,6,9,10,11,14}
\lvs{1,4,5,6,9,10,11,14}
\end{scope}
\end{tikzpicture}
\caption{The generators $\ol f_{i,j}\in\Sing(\I_n)$ for $i<j$ (left) and $i>j$ (right).}
\label{fig:f}
\end{center}
\end{figure}
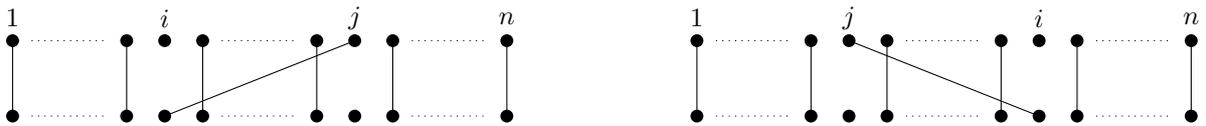

Let $\RS$ be the following set of relations, where $i,j,k,l\in\bn$ are distinct in each:
\begin{align}
\label{SI1} f_{i,j}f_{j,i}f_{i,j} &= f_{i,j}, \\
\label{SI2} f_{i,j}^3=f_{i,j}^2 &= f_{j,i}^2, \\
\label{SI3} f_{i,j}f_{k,l} &= f_{k,l}f_{i,j}, \\
\label{SI4} f_{i,j}f_{j,i} &= f_{i,k}f_{k,i}, \\
\label{SI7} f_{i,j}f_{i,k}=f_{j,k}f_{i,j} &= f_{i,k}f_{j,k}, \\
\label{SI5} f_{k,i}f_{i,j}f_{j,k} &= f_{k,j}f_{j,i}f_{i,k}, \\
\label{SI6} f_{k,i}f_{i,j}f_{j,k}f_{k,l} &= f_{k,l}f_{l,i}f_{i,j}f_{j,l}.
\end{align}
The next result is \cite[Theorem 2.1]{East2015}, but we note that the $f_{i,j}$ were called $e_{ij}$ in \cite{East2015}:

\begin{theorem}\label{thm:SI}
For $n\geq2$, the semigroup $\Sing(\I_n)$ has presentation $\Spres{\XS}{\RS}$ via~$\pS$.  \epfres
\end{theorem}

\subsection[Presentation for $\Sing(M_0^n)$]{\boldmath Presentation for $\Sing(M_0^n)$}

In this subsection we aim to give a presentation for the semigroup
\[
\Sing(M_0^n) = M_0^n\sm M^n
\]
of all $n$-tuples over $M_0$ with at least one zero entry.  The notation $\Sing(M_0^n)$ is for convenience only, and is not meant to imply that $M^n$ is the group of units of $M_0^n$ (which is only the case when~$M$ is a group).  We continue to fix the monoid presentation $\Mpres{X_M}{R_M}$ for $M$ via $\phi_M:X_M^*\to M$, and we write $\ol w=w\phi_M$ for $w\in X_M^*$. Since we are incorporating this \emph{monoid} presentation into a \emph{semigroup} presentation for $\Sing(M_0^n)$, we will have to tread carefully when considering empty words.

We begin by defining the alphabet 
\begin{equation}\label{eq:XSM}
\XSM = \{e_1,\ldots,e_n\} \cup \set{x^{(i;j)}}{x\in X_M,\ i,j\in\bn,\ i\not=j},
\end{equation}
and the morphism
\begin{equation}\label{eq:pSM}
\pSM : \XSM^+\to\Sing(M_0^n) : e_i\mt\ol e_i , \ x^{(i;j)}\mt\ol x^{(i;j)},
\end{equation}
which is defined diagrammatically in Figure~\ref{fig:ex}.  Here as usual we identify $\Sing(M_0^n)$ with a subsemigroup of ${M\wr\Sing(\I_n)}$, and we recall the usual convention about unspecified upper vertex labels.  We note that $\ol x^{(i;j)} = \ol x^{(i)}\ol e_j = \ol e_j\ol x^{(i)}$ in the notation of Subsection~\ref{subsect:MIn}, but we also note that $\ol x^{(i)}\not\in\Sing(M_0^n)$.  It is worth noting that
\[
\supp(\ol e_i) = \{i\}\com \AND \supp(\ol x^{(i;j)}) = \{j\}\com
\]
for appropriate $i,j,x$, again using the $A\com=\bn\sm A$ notation.
As usual, we write $\ol w=w\pSM$ for $w\in\XSM^+$.

\begin{figure}[ht]
\begin{center}
\begin{tikzpicture}[scale=.5]
\begin{scope}[shift={(0,0)}]	
\udotted14
\ddotted14
\udotted69
\ddotted69
\udotted{11}{14}
\ddotted{11}{14}
\stline11
\stline44
\stline66
\stline99
\stline55
\stline{11}{11}
\stline{14}{14}
\node()at(1,2.6){\small$1$};
\node()at(5,2.9){\small$i$};
\node()at(10,2.6){\small$j$};
\node()at(14,2.6){\small$n$};
\uvs{1,4,5,6,9,10,11,14}
\lvs{1,4,5,6,9,10,11,14}
\node[circle,draw=black, fill=white, inner sep = 0.06cm] () at (5,2){\footnotesize$\ol x$};
\end{scope}
\begin{scope}[shift={(19,0)}]	
\udotted14
\ddotted14
\udotted69
\ddotted69
\udotted{11}{14}
\ddotted{11}{14}
\stline11
\stline44
\stline66
\stline99
\stline{10}{10}
\stline{11}{11}
\stline{14}{14}
\node()at(1,2.6){\small$1$};
\node()at(5,2.6){\small$j$};
\node()at(10,2.9){\small$i$};
\node()at(14,2.6){\small$n$};
\uvs{1,4,5,6,9,10,11,14}
\lvs{1,4,5,6,9,10,11,14}
\node[circle,draw=black, fill=white, inner sep = 0.06cm] () at (10,2){\footnotesize$\ol x$};
\end{scope}
\begin{scope}[shift={(12,4)}]	
\uvs{1,4,5,6,9}
\lvs{1,4,5,6,9}
\udotted14
\ddotted14
\udotted69
\ddotted69
\stline11
\stline44
\stline66
\stline{9}{9}
\node()at(1,2.6){\small$1$};
\node()at(5,2.6){\small$i$};
\node()at(9,2.6){\small$n$};
\end{scope}
\end{tikzpicture}
\caption{Generators of $\Sing(M_0^n)$: $\ol e_i$ (top), and $\ol x^{(i;j)}$ for $i<j$ (left) and $i>j$ (right).}
\label{fig:ex}
\end{center}
\end{figure}
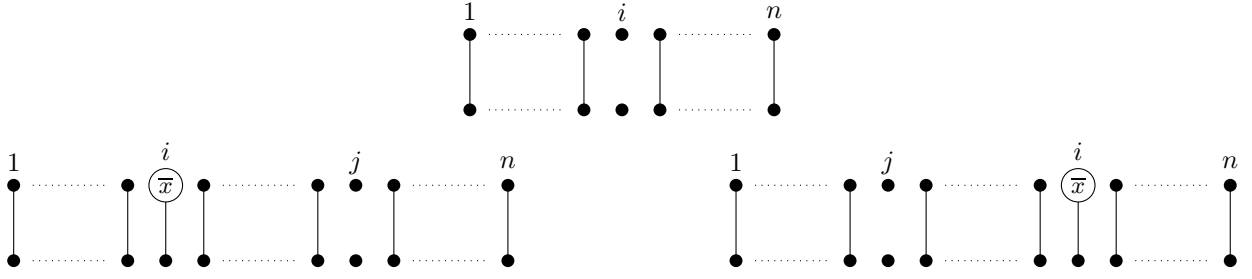

For a (possibly empty) word $w=x_1\cdots x_k\in X_M^*$, and for distinct $i,j\in\bn$, we define
\[
w^{(i;j)} = \begin{cases}
x_1^{(i;j)}\cdots x_k^{(i;j)} &\text{if $k\geq1$,}\\
e_j &\text{if $k=0$.}
\end{cases}
\]
Note that $w^{(i;j)}$ is always a non-empty word over $\XSM$.  We now define $\RSM$ to be the following set of relations, where $x,y\in X_M$ and $i,j,k,l\in\bn$ range over all meaningful values, subject to any stated constraints:
\begin{align}
\label{M1} e_i^2 &= e_i, \\
\label{M2} e_ie_j &= e_je_i, \\
\label{M3} u^{(i;j)} &= v^{(i;j)} &&\text{for $(u,v)\in R_M$,}\\
\label{M4} x^{(i;j)}e_k &= e_kx^{(i;j)} = e_jx^{(i;k)} &&\text{if $k\not=i,j$,} \\
\label{M5} x^{(i;j)}e_j &= e_jx^{(i;j)} = x^{(i;j)} , \\
\label{M6} x^{(i;j)}e_i &= e_ix^{(i;j)} = e_ie_j , \\
\label{M7} x^{(i;k)}y^{(j;k)} &= y^{(j;k)}x^{(i;k)} &&\text{if $i\not=j$.}
\end{align}
For the following proof it is convenient to introduce a piece of notation.  For $a\in M$, and for distinct $i,j\in\bn$, we define the tuple
\[
a^{(i;j)} = (b_1,\ldots,b_n) \WHERE b_k = \begin{cases}
a &\text{if $k=i$,}\\
0 &\text{if $k=j$,}\\
1 &\text{otherwise.}
\end{cases}
\]
So $a^{(i;j)}\in\Sing(M_0^n)$ for all such $a,i,j$.  When $a=\ol x$ for some $x\in X_M$, this agrees with our earlier definition of $\ol x^{(i;j)}$; see Figure \ref{fig:ex}.

\begin{theorem}\label{thm:M}
For $n\geq2$, the semigroup $\Sing(M_0^n)=M_0^n\sm M^n$ has presentation
\[
\Spres{\XSM}{\RSM}
\]
via $\pSM$.
\end{theorem}

\pf
To show that $\pSM$ is surjective, let $\ba=(a_1,\ldots,a_n)\in\Sing(M_0^n)$; we must show that $\ba=\ol w$ for some $w\in\XSM^+$.  Re-naming the elements of $\bn$ if necessary, we may assume without loss of generality that $\supp(\ba)=\{1,\ldots,q\}$ for some $0\leq q<n$.  So in fact, $\ba=(a_1,\ldots,a_q,0,\ldots,0)$ with $a_1,\ldots,a_q\in M$.  For $1\leq i\leq q$, we have $a_i=w_i\phi_M$ for some $w_i\in X_M^*$, and it follows that
\[
\ba = \ol e_{q+1}\cdots\ol e_n \cdot a_1^{(1;n)}\cdots a_q^{(q;n)} = \ol{e_{q+1}\cdots e_n\cdot w_1^{(1;n)}\cdots w_q^{(q;n)}},
\]
as required.  (When $q=0$, the above says that $\ba=\ol{e_1\cdots e_n}$.  If $w_i=\io$ for some $1\leq i\leq q$, then $a_i=\ol w_i = 1$, and $a_i^{(i;n)} = 1^{(i;n)} = \ol e_n = \ol{\io^{(i;n)}} = \ol{w_i^{(i;n)}}$.)

Next, one can easily check that the relations from $\RSM$ are preserved by $\pSM$; see Figure \ref{fig:M} for an example calculation for relation \eqref{M4}.  It follows that ${\RSM^\sharp\sub\ker(\pSM)}$.  The rest of the proof involves establishing the reverse inclusion.  The bulk of the work goes into proving the following technical lemma.  For the rest of the proof we write ${\sim}$ for $\RSM^\sharp$.

\begin{lemma}\label{lem:M}
If $w\in\XSM^+$ is such that $\supp(\ol w)=\{1,\ldots,q\}$ for some $0\leq q<n$, then
\[
w\sim e_{q+1}\cdots e_n\cdot w_1^{(1;n)}\cdots w_q^{(q;n)} \qquad\text{for some $w_1,\ldots,w_q\in X_M^*$.}
\]
\end{lemma}

\pf
Before we begin the proof of the lemma, we note that while some of the $w_i$ ($1\leq i\leq q$) could be empty, the $w_i^{(i;n)}$ never are.  Indeed, if $w_i=\io$ for some $i$, then $w_i^{(i;n)} = \io^{(i;n)} = e_n$.

Beginning the proof now, we first observe that (either part of) Lemma \ref{lem:XY} and relations \mbox{\eqref{M4}--\eqref{M6}} give
\[
w\sim uv \qquad\text{for some $u\in\{e_1,\ldots,e_n\}^*$ and $v\in\set{x^{(i;j)}}{x\in X_M,\ i,j\in\bn,\ i\not=j}^*$.}
\]
Since $w$ is non-empty, certainly $u$ and $v$ are not both empty, and we write
\[
u = e_{i_1}\cdots e_{i_s} \AND v = x_1^{(j_1;k_1)} \cdots x_t^{(j_t;k_t)}.
\]
Because of \eqref{M1} and \eqref{M2} we can assume that $i_1<\cdots<i_s$.  Now,
\begin{align*}
\{1,\ldots,q\} = \supp(\ol w) &= \supp(\ol e_{i_1}) \cap\cdots\cap \supp(\ol e_{i_s}) \cap \supp(\ol x_1^{(j_1;k_1)})\cap \cdots \cap\supp(\ol x_t^{(j_t;k_t)})\\
&=\{i_1\}\com\cap\cdots\cap\{i_s\}\com\cap\{k_1\}\com\cap\cdots\cap\{k_t\}\com
= \{i_1,\ldots,i_s,k_1,\ldots,k_t\}\com,
\end{align*}
so it follows that
\begin{equation}\label{eq:qn}
\{i_1,\ldots,i_s,k_1,\ldots,k_t\} = \{q+1,\ldots,n\}.  
\end{equation}

Next we claim that
\begin{equation}\label{eq:u'v}
w\sim u' v \WHERE u' = e_{q+1}\cdots e_n.
\end{equation}
Indeed, this is obvious if $\{i_1,\ldots,i_s\}=\{q+1,\ldots,n\}$, so suppose instead that there exists some $l\in\{q+1,\ldots,n\}\sm\{i_1,\ldots,i_s\}$.  By \eqref{eq:qn} we have $l=k_m$ for some $1\leq m\leq t$.  For convenience we write
\begin{equation}\label{eq:v1v2}
v_1 = x_1^{(j_1;k_1)} \cdots x_{m-1}^{(j_{m-1};k_{m-1})} \AND v_2 = x_{m+1}^{(j_{m+1};k_{m+1})} \cdots x_t^{(j_t;k_t)}.
\end{equation}
Then by \eqref{M4}--\eqref{M6} we have
\[
v = v_1 \cdot x_m^{(j_m;l)} \cdot v_2 \sim v_1 \cdot e_lx_m^{(j_m;l)} \cdot v_2 \sim e_l\cdot v_1 \cdot x_m^{(j_m;l)} \cdot v_2 = e_l \cdot v.
\]
It follows that $w\sim uv \sim ue_l\cdot v$.  Continuing in this way, we have
\[
w\sim ue_{l_1}\cdots e_{l_r}\cdot v \WHERE \{q+1,\ldots,n\}\sm\{i_1,\ldots,i_s\} = \{l_1,\ldots,l_r\}.
\]
Since
\[
ue_{l_1}\cdots e_{l_r} = e_{i_1}\cdots e_{i_s}\cdot e_{l_1}\cdots e_{l_r} \sim e_{q+1}\cdots e_n = u',
\]
by \eqref{M2}, this completes the proof of \eqref{eq:u'v}.  In what follows, we continue to write ${u' = e_{q+1}\cdots e_n}$.

Next we claim that
\begin{equation}\label{eq:u'v'}
w\sim u' v'  \qquad\text{for some $v'\in\set{x^{(i;n)}}{x\in X_M,\ 1\leq i<n}^*$.}
\end{equation}
This follows immediately from \eqref{eq:u'v} if $k_m=n$ for all $1\leq m\leq t$ (in which case we take $v'=v$), so suppose instead that $k_m\not=n$ for some such $m$.  We prove \eqref{eq:u'v'} by showing that we can use the relations to either remove the letter $x_m^{(j_m;k_m)}$ from $v$ or else replace it by $x_m^{(j_m;n)}$.
Again for convenience we let~$v_1$ and~$v_2$ be as in~\eqref{eq:v1v2}.  Continuing from \eqref{eq:u'v}, and using relations \eqref{M1} and \eqref{M4}--\eqref{M6}, we first note that
\begin{equation}\label{eq:wuvexv}
w \sim u'v \sim u'e_n \cdot v_1\cdot x_m^{(j_m;k_m)}\cdot v_2 \sim u' \cdot v_1\cdot e_nx_m^{(j_m;k_m)}\cdot v_2.
\end{equation}
From this point we must consider separate cases for $j_m=n$ and $j_m\not=n$.  For both cases, we use relations \eqref{M1}, \eqref{M2} and \eqref{M4}--\eqref{M6}, and we keep $q+1\leq k_m< n$ in mind; cf.~\eqref{eq:qn}.  If $j_m=n$, then continuing from \eqref{eq:wuvexv}, we have
\[
w \sim u' \cdot v_1\cdot e_ne_{k_m}\cdot v_2 
\sim u'e_ne_{k_m} \cdot v_1\cdot v_2 
\sim u' \cdot v_1\cdot v_2,
\]
so we have deleted $x_m^{(j_m;k_m)}$ in this case.  On the other hand, if $j_m\not=n$, then continuing from~\eqref{eq:wuvexv}, we have
\[
w \sim u'\cdot v_1\cdot e_nx_m^{(j_m;k_m)}\cdot v_2 
\sim u'\cdot v_1\cdot e_{k_m}x_m^{(j_m;n)}\cdot v_2 
\sim u'e_{k_m}\cdot v_1\cdot x_m^{(j_m;n)}\cdot v_2 
\sim u'\cdot v_1\cdot x_m^{(j_m;n)}\cdot v_2,
\]
so we have replaced $x_m^{(j_m;k_m)}$ by $x_m^{(j_m;n)}$ in this case.  As noted above, this completes the proof of~\eqref{eq:u'v'}.

Next we claim that
\begin{equation}\label{eq:u'v''}
w\sim u' v''  \qquad\text{for some $v''\in\set{x^{(i;n)}}{x\in X_M,\ 1\leq i\leq q}^*$.}
\end{equation}
To prove this, we continue from \eqref{eq:u'v'}, and we show that any letter $x^{(i;n)}$ appearing in $v'$ with $q<i<n$ can be deleted (using the relations).  To do so, suppose some such letter exists, and write $v'=v_1'\cdot x^{(i;n)}\cdot v_2'$.  Then using relations \eqref{M1}, \eqref{M2} and \eqref{M4}--\eqref{M6}, we have
\[
w\sim u'v' \sim u'e_i \cdot v_1'\cdot x^{(i;n)}\cdot v_2' \sim u' \cdot v_1'\cdot e_ix^{(i;n)}\cdot v_2' \sim u' \cdot v_1'\cdot e_ie_n\cdot v_2' \sim u'e_ie_n \cdot v_1'\cdot  v_2' \sim u' \cdot v_1'\cdot  v_2',
\]
as required.

Next, we use \eqref{M7} to `unshuffle' $v''$ to obtain
\begin{equation}\label{eq:wuvv}
w \sim u'v'' \sim u'\cdot v_1\cdots v_q \qquad\text{where $v_i\in\set{x^{(i;n)}}{x\in X_M}^*$ for each $1\leq i\leq q$.}
\end{equation}
If some $v_i$ is non-empty, then of course it has the form $w_i^{(i;n)}$ for some $w_i\in X_M^+$.  In particular, if every $v_i$ is non-empty, then we have proved the lemma.  So now suppose some $v_i=\io$ is empty.  Then by \eqref{M1} and \eqref{M5} we have
\[
w \sim u' \cdot v_1\cdots v_{i-1}\cdot v_{i+1}\cdots v_q 
\sim u'e_n \cdot v_1\cdots v_{i-1}\cdot v_{i+1}\cdots v_q 
\sim u' \cdot v_1\cdots v_{i-1}\cdot e_n\cdot v_{i+1}\cdots v_q ,
\]
so we can replace $v_i=\io$ by $e_n=\io^{(i;n)}$ in \eqref{eq:wuvv}.  After doing this for every $i$ for which $v_i=\io$, the lemma has been proved.
\epf

Before returning to the main proof, we also record the following simple lemma.  For the proof, we observe that for any $u,v\in X_M^*$ and distinct $i,j\in\bn$, we have $(uv)^{(i;j)}\sim u^{(i;j)}v^{(i;j)}$.  Indeed, this follows immediately from the definitions if $u$ and $v$ are both non-empty, from \eqref{M1} if both are empty, or from \eqref{M5} if only one is empty.

\begin{lemma}\label{lem:uv}
If $(u,v)\in R_M^\sharp$, and if $i,j\in\bn$ are distinct, then $u^{(i;j)}\sim v^{(i;j)}$.
\end{lemma}

\pf
It suffices to consider the case that $u$ and $v$ differ by a single application of a relation from~$R_M$.  So we assume that
\[
u = su't \AND v = sv't \qquad\text{for some $s,t\in X_M^*$ and $(u',v')\in R_M$.}
\]
But then combining the observation before the lemma with \eqref{M3}, we have
\[
u^{(i;j)} \sim s^{(i;j)}(u')^{(i;j)}t^{(i;j)} \sim s^{(i;j)}(v')^{(i;j)}t^{(i;j)} \sim v^{(i;j)},
\]
as required.
\epf

Returning now to the proof of the theorem, suppose $(u,v)\in\ker(\pSM)$.  So ${u,v\in\XSM^+}$ and $\ol u=\ol v$; we must show that $u\sim v$.  Re-naming the elements of $\bn$ if necessary, we may assume for convenience that $\supp(\ol u)=\supp(\ol v)=\{1,\ldots,q\}$ for some $0\leq q<n$.  By Lemma \ref{lem:M}, we have
\[
u\sim e_{q+1}\cdots e_n\cdot u_1^{(1;n)}\cdots u_q^{(q;n)} \AND v\sim e_{q+1}\cdots e_n\cdot v_1^{(1;n)}\cdots v_q^{(q;n)} \qquad\text{for some $u_i,v_i\in X_M^*$.}
\]
But then
\[
(\ol u_1,\ldots,\ol u_q,0,\ldots,0) = \ol u = \ol v = (\ol v_1,\ldots,\ol v_q,0,\ldots,0),
\]
and so $\ol u_i=\ol v_i$ for all $1\leq i\leq q$.  Since $M$ has presentation $\Mpres{X_M}{R_M}$, it follows that $(u_i,v_i)\in R_M^\sharp$ for all $i$.  Lemma \ref{lem:uv} then gives $u_i^{(i;n)}\sim v_i^{(i;n)}$ for all $i$.  Putting everything together, we finally deduce that
\[
u \sim e_{q+1}\cdots e_n\cdot u_1^{(1;n)}\cdots u_q^{(q;n)} \sim e_{q+1}\cdots e_n\cdot v_1^{(1;n)}\cdots v_q^{(q;n)} \sim v,
\]
and this completes the proof.
\epf

\begin{figure}[ht]
\begin{center}
\begin{tikzpicture}[scale=.7]
\begin{scope}[shift={(0,0)}]	
\stline11
\stline33
\uvs{1,2,3}
\lvs{1,2,3}
\node () at (4.5,0) {$=$};
\node[circle,draw=black, fill=white, inner sep = 0.06cm] () at (1,2){\footnotesize$\ol x$};
\node()at(1,2.6){\small$i$};
\node()at(2,2.6){\small$j$};
\node()at(3,2.6){\small$k$};
\end{scope}
\begin{scope}[shift={(0,-2)}]	
\stline11
\stline22
\uvs{1,2,3}
\lvs{1,2,3}
\end{scope}
\begin{scope}[shift={(5,0)}]	
\stline11
\stline22
\uvs{1,2,3}
\lvs{1,2,3}
\node () at (4.5,0) {$=$};
\node()at(1,2.6){\small$i$};
\node()at(2,2.6){\small$j$};
\node()at(3,2.6){\small$k$};
\end{scope}
\begin{scope}[shift={(5,-2)}]	
\stline11
\stline33
\uvs{1,2,3}
\lvs{1,2,3}
\node[circle,draw=black, fill=white, inner sep = 0.06cm] () at (1,2){\footnotesize$\ol x$};
\end{scope}
\begin{scope}[shift={(10,0)}]	
\stline11
\stline33
\uvs{1,2,3}
\lvs{1,2,3}
\node()at(1,2.6){\small$i$};
\node()at(2,2.6){\small$j$};
\node()at(3,2.6){\small$k$};
\end{scope}
\begin{scope}[shift={(10,-2)}]	
\stline11
\stline22
\uvs{1,2,3}
\lvs{1,2,3}
\node[circle,draw=black, fill=white, inner sep = 0.06cm] () at (1,2){\footnotesize$\ol x$};
\end{scope}
\end{tikzpicture}
\caption{Relation \eqref{M4}, picturing only strings $i,j,k$.  (Note that no ordering on $i,j,k$ is implied.)}
\label{fig:M}
\end{center}
\end{figure}
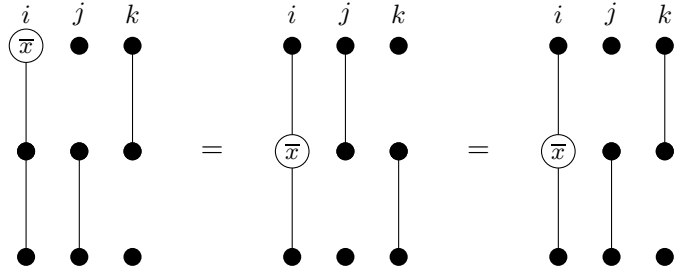

\subsection[Presentation for $M\wr\Sing(\I_n)$]{\boldmath Presentation for $M\wr\Sing(\I_n)$}

We are now ready to give our presentation for the singular wreath product $M\wr\Sing(\I_n)$.  We begin by defining the alphabet
\[
\XMS = \XS \cup \XSM , 
\]
where $\XS$ and $\XSM$ are as in \eqref{eq:XS} and \eqref{eq:XSM}.  We define the morphism
\[
\pMS : \XMS^+\to M\wr\Sing(\I_n)
\]
to be the extension of the morphisms
\[
\pS : \XS^+\to \Sing(\I_n) \AND \pSM : \XSM^+\to \Sing(M_0^n)
\]
from \eqref{eq:pS} and \eqref{eq:pSM}, where as usual we identify $\Sing(\I_n)$ and $\Sing(M_0^n)$ with subsemigroups of $M\wr\Sing(\I_n)$.

Finally, we let $\RMS$ be the set of relations over $\XMS$ consisting of ${\RS\cup\RSM}$ (i.e., \eqref{SI1}--\eqref{SI6} and \eqref{M1}--\eqref{M7}), together with the following, where $x\in X_M$ and $i,j,k,l\in\bn$ range over all meaningful values in each relation:
\begin{align}
\label{MS1} f_{i,j}f_{j,i} &= e_i,\\
\label{MS2} f_{i,j}x^{(i;j)} &= x^{(j;i)}f_{i,j}, \\
\label{MS3} f_{i,j}x^{(i;k)} &= x^{(j;k)}f_{i,j}, \\
\label{MS3'} f_{i,j}x^{(k;i)} &= x^{(k;i)}e_j, \\
\label{MS4} f_{i,j}x^{(j;k)} &= f_{i,j}f_{k,j}, \\
\label{MS4'} f_{i,j}x^{(k;j)} &= x^{(k;i)}f_{i,j}, \\
\label{MS5} f_{i,j} x^{(k;l)} &= x^{(k;l)}f_{i,j} &&\hspace{-2cm}\text{if $\{i,j\}\cap\{k,l\}=\es$.}
\end{align}
Note for example that the `all meaningful values' assumption implies that $i,j,k$ are distinct in \eqref{MS3}.  On the other hand, \eqref{MS4} includes the case in which $i,j,k$ are distinct, as well as the case that $k=i\not=j$.

\begin{theorem}\label{thm:MSIn}
For $n\geq2$, the singular wreath product $M\wr\Sing(\I_n)$ has presentation
\[
\Spres{\XMS}{\RMS}
\]
via $\pMS$.
\end{theorem}

\pf
Surjectivity of $\pMS$ follows quickly from the surjectivity of $\pS$ and $\pSM$, given that $(\ba,\al)=\ba\cdot\al$ for any $(\ba,\al)\in M\wr\Sing(\I_n)$; cf.~Remark \ref{rem:aabb}.  One can also easily check that $\pMS$ preserves the relations; for example, Figure \ref{fig:MS} does this for \eqref{MS4}.  This shows that ${\RMS^\sharp\sub\ker(\pMS)}$.  Thus, as ever, the bulk of the proof goes into showing the reverse inclusion.  For the rest of the proof, we write ${\sim}$ for $\RMS^\sharp$.

\begin{lemma}\label{lem:w'w''}
For any $w\in\XMS^+$, we have
\[
w\sim w'w'' \qquad\text{for some $w'\in\XSM^+$ and $w''\in\XS^+$ with $\supp(\ol w')=\dom(\ol w'')$.}
\]
\end{lemma}

\pf
We first use Lemma \ref{lem:XY}\ref{XY1} with $X=\XSM$ and $Y=\XS$ to deduce that
\begin{equation}\label{eq:wuv}
w \sim uv \qquad\text{for some $u\in\XSM^*$ and $v\in\XS^*$.}
\end{equation}
To show that the lemma applies, we need to check that for each $x\in\XSM$ and $y\in\XS$, we have $yx\sim uv$ for some $u\in\XSM^*$ and $v\in\XS^*$ with $\ell(u)\leq1$.  Now, $y$ of course has the form $y=f_{i,j}$.  The case that $x=e_k$ follows from \eqref{MS1}, as then
\[
yx = f_{i,j}e_k \sim f_{i,j}f_{k,j}f_{j,k},
\]
and we take $u=\io$ and $v=f_{i,j}f_{k,j}f_{j,k}$.  For the case that $x=x^{(k;l)}$, we just need to observe that \eqref{MS2}--\eqref{MS5} contains a relation $(yx,uv)$ of the desired form in each case.  We now fix $u$ and $v$ as in \eqref{eq:wuv}.  Next we claim that
\begin{equation}\label{eq:wu'v}
w \sim w'v \qquad\text{for some $w'\in\XSM^+$ with $\supp(\ol w')\sub\dom(\ol v)$.}
\end{equation}
Indeed, there is nothing to show if $\dom(\ol v)=\bn$ (i.e., if $v=\io$), so suppose instead that ${\dom(\ol v)\not=\bn}$, and fix some $i\in\dom(\ol v)\com$.  Since $n\geq2$, we may also fix some $j\in\{i\}\com$.  Since $i\not\in\dom(\ol v)$, and since $\ol f_{i,j}\ol f_{j,i}=\id_{\{i\}\com}$, we have $\ol v=\ol f_{i,j}\ol f_{j,i}\ol v$.  It follows from Theorem \ref{thm:SI} that $v\sim f_{i,j}f_{j,i}v$.  We then combine this with \eqref{MS1} to calculate
\[
w \sim uv \sim uf_{i,j}f_{j,i}\cdot v \sim ue_i \cdot v.
\]
Continuing in this way, and writing $\dom(\ol v)\com = \{i_1,\ldots,i_k\}$, we have
\[
w \sim ue_{i_1}\cdots e_{i_k}\cdot v  .
\]
We see then that \eqref{eq:wu'v} holds with $w'=ue_{i_1}\cdots e_{i_k}$, as
\begin{align*}
\supp(\ol w') &= \supp(\ol u)\cap\supp(\ol e_{i_1})\cap\cdots\cap\supp(\ol e_{i_k})\\
&= \supp(\ol u) \cap \{i_1\}\com \cap\cdots\cap\{i_k\}\com \\
&= \supp(\ol u) \cap \{i_1,\ldots,i_k\}\com \sub \{i_1,\ldots,i_k\}\com = \dom(\ol v).
\end{align*}

We now fix $w'$ and $v$ as in \eqref{eq:wu'v}.  Also fix some $i\in\supp(\ol w')\com$, and some $j\in\{i\}\com$.  Since $i\not\in\supp(\ol w')$, we have $\ol w'=\ol w'\ol e_i$, and so $w'\sim w'e_i$ by Theorem \ref{thm:M}.  Combining this with~\eqref{MS1}, it follows that
\[
w \sim w'v \sim w'\cdot e_iv \sim w' \cdot f_{i,j}f_{j,i}v.
\]
Continuing in this way, writing $\supp(\ol w')\com=\{i_1,\ldots,i_k\}$, and fixing some $j_s\in\{i_s\}\com$ for each $1\leq s\leq k$, we have
\[
w \sim w' \cdot (f_{i_1,j_1}f_{j_1,i_1})\cdots(f_{i_k,j_k}f_{j_k,i_k})\cdot v.
\]
We see then that the lemma holds with $w''=(f_{i_1,j_1}f_{j_1,i_1})\cdots(f_{i_k,j_k}f_{j_k,i_k})\cdot v$.  Indeed, we first note that 
\[
\ol{(f_{i_1,j_1}f_{j_1,i_1})\cdots(f_{i_k,j_k}f_{j_k,i_k})} = \id_{\{i_1\}\com} \cdots\id_{\{i_k\}\com}  = \id_{\{i_1,\ldots,i_k\}\com}  = \id_{\supp(\ol w')}.
\]
But then
\[
\dom(\ol w'') = \dom(\id_{\supp(\ol w')} \cdot \ol v) = \dom(\id_{\supp(\ol w')} ) \cap \dom(\ol v) = \supp(\ol w') \cap \dom(\ol v) = \supp(\ol w'),
\]
where the last equality follows from the fact that $\supp(\ol w')\sub\dom(\ol v)$; cf.~\eqref{eq:wu'v}.
\epf

Returning to the proof of the theorem, let $(u,v)\in\ker(\pMS)$, so that $u,v\in\XMS^+$ and $\ol u=\ol v$; we must show that $u\sim v$.  By Lemma \ref{lem:w'w''}, we have
\[
u\sim u'u'' \AND v\sim v'v'' 
\]
for some $u',v'\in\XSM^+$ and $u'',v''\in\XS^+$ with $\supp(\ol u')=\dom(\ol u'')$ and ${\supp(\ol v')=\dom(\ol v'')}$.  Keeping Remark \ref{rem:aabb} in mind, it follows that
\[
(\ol u',\ol u'') = \ol u'\cdot\ol u'' = \ol u = \ol v = \ol v'\cdot\ol v'' = (\ol v',\ol v''),
\]
and so $\ol u'=\ol v'$ (in $\Sing(M_0^n)$) and $\ol u''=\ol v''$ (in $\Sing(\I_n)$).  It then follows from Theorems \ref{thm:SI} and~\ref{thm:M} that $u'\sim v'$ and $u''\sim v''$, and so $u \sim u'u'' \sim v'v'' \sim v$, completing the proof of the theorem.
\epf

\begin{figure}[ht]
\begin{center}
\begin{tikzpicture}[scale=.7]
\begin{scope}[shift={(0,0)}]	
\stline21
\uvs{1,2}
\lvs{1,2}
\node () at (3.5,0) {$=$};
\node()at(1,2.6){\small$k=i\phantom{{}=k}$};
\node()at(2,2.6){\small$j$};
\end{scope}
\begin{scope}[shift={(0,-2)}]	
\stline22
\uvs{1,2}
\lvs{1,2}
\node[circle,draw=black, fill=white, inner sep = 0.06cm] () at (2,2){\footnotesize$\ol x$};
\end{scope}
\begin{scope}[shift={(4,0)}]	
\stline21
\uvs{1,2}
\lvs{1,2}
\node()at(1,2.6){\small$k=i\phantom{{}=k}$};
\node()at(2,2.6){\small$j$};
\end{scope}
\begin{scope}[shift={(4,-2)}]	
\stline21
\uvs{1,2}
\lvs{1,2}
\end{scope}
\begin{scope}[shift={(10,0)}]	
\stline21
\stline33
\uvs{1,2,3}
\lvs{1,2,3}
\node () at (4.5,0) {$=$};
\node()at(1,2.6){\small$i$};
\node()at(2,2.6){\small$j$};
\node()at(3,2.6){\small$k$};
\end{scope}
\begin{scope}[shift={(10,-2)}]	
\stline11
\stline22
\uvs{1,2,3}
\lvs{1,2,3}
\node[circle,draw=black, fill=white, inner sep = 0.06cm] () at (2,2){\footnotesize$\ol x$};
\end{scope}
\begin{scope}[shift={(15,0)}]	
\stline21
\stline33
\uvs{1,2,3}
\lvs{1,2,3}
\node()at(1,2.6){\small$i$};
\node()at(2,2.6){\small$j$};
\node()at(3,2.6){\small$k$};
\end{scope}
\begin{scope}[shift={(15,-2)}]	
\stline11
\stline23
\uvs{1,2,3}
\lvs{1,2,3}
\end{scope}
\end{tikzpicture}
\caption{Relation \eqref{MS4}, picturing only strings $i,j,k$, and considering the cases $i=k$ (left) and $i\not=k$ (right).  (Note that no ordering on $i,j,k$ is implied.)}
\label{fig:MS}
\end{center}
\end{figure}
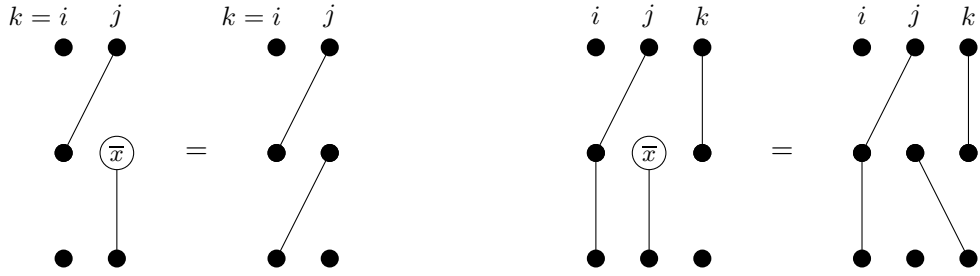

In light of relation \eqref{MS1}, one could remove the generators $e_i$ ($i\in\bn$) from the presentation $\Spres{\XMS}{\RMS}$, replacing each occurrence of $e_i$ in the relations by some fixed word of the form $f_{i,j}f_{j,i}$ ($j\in\{i\}\com$).  We note however, that the resulting set of relations would not be quite as symmetrical as $\RMS$.

\footnotesize
\def\bibspacing{-1.1pt}
\bibliography{bibRevised}
\bibliographystyle{abbrv}

\end{document}